# Preconditioning of a Generalized Forward-Backward Splitting and Application to Optimization on Graphs

**Hugo Raguet**
hugo.raguet@gmail.com

**Loïc Landrieu**
ENPC, *IMAGINE/LIGM*
INRIA, *SIERRA*
landrieu.loic@gmail.com

**Abstract**

We present a preconditioning of a generalized forward-backward splitting algorithm for finding a zero of a sum of maximally monotone operators $\sum_{i=1}^{n} A_i + B$ with $B$ cocoercive, involving only the computation of $B$ and of the resolvent of each $A_i$ separately. This allows in particular to minimize functionals of the form $\sum_{i=1}^{n} g_i + f$ with $f$ smooth, using only the gradient of $f$ and the proximity operator of each $g_i$ separately. By adapting the underlying metric, such preconditioning can serve two practical purposes: first, it might accelerate the convergence, or second, it might simplify the computation of the resolvent of $A_i$ for some $i$. In addition, in many cases of interest, our preconditioning strategy allows the economy of storage and computation concerning some auxiliary variables. In particular, we show how this approach can handle large-scale, nonsmooth, convex optimization problems structured on graphs, which arises in many image processing or learning applications, and that it compares favorably to alternatives in the literature.

## 1 Introduction

The generalization of the forward-backward splitting algorithm which we presented recently (Raguet *et al.*, 2013) aims at solving, over a real Hilbert space $\mathcal{H}$, monotone inclusion problems of the form (see § 2.1 for notations)

$$\text{find } x \in \operatorname{zer}\left\{\sum_{i=1}^{n} A_i + B\right\}, \tag{1.1}$$

by making only use of the *resolvent* of each set-valued operator $A_i\colon \mathcal{H} \to 2^{\mathcal{H}}$, supposed maximally monotone, and of the explicit application of $B\colon \mathcal{H} \to \mathcal{H}$, supposed *cocoercive*. Namely, noting $x$ the main iterate and introducing $n$ auxiliary variables $(z_i)_{1 \le i \le n} \in \mathcal{H}^n$, iteration $k \in \mathbb{N}$ of the algorithm consists in

$$\begin{array}{l} \textbf{for } i \in \{1, \dots, n\} \textbf{ do} \\ \quad \left\lfloor\; z_i \leftarrow z_i + \rho_k\Big(J_{\frac{\gamma}{w_i} A_i}(2x - z_i - \gamma B x) - x\Big)\,, \right. \\ x \leftarrow \sum_{i=1}^{n} w_i z_i\,, \end{array} \tag{1.2}$$

where for each $i$, $J_{\frac{\gamma}{w_i}A_i} \stackrel{\text{def}}{=} \left(\mathrm{Id} + \frac{\gamma}{w_i}A_i\right)^{-1}$ denotes the resolvent of $\frac{\gamma}{w_i}A_i$. The weights $(w_i)_{1\leq i\leq n}$ sum up to unity and tune the influence of each operator in the splitting; the step size $\gamma$ is constant along iterations, while each $\rho_k$ is a relaxation parameter.

Interestingly, the case $n \stackrel{\text{set}}{=} 1$ reduces to the relaxed *forward-backward* splitting algorithm (Combettes, 2004, Section 6); in this case, convergence can be ensured with step size $\gamma$ varying along iterations. On the other hand, the case $B \stackrel{\text{set}}{=} 0$ retrieves the relaxed *Douglas–Rachford* splitting algorithm (Combettes, 2004, Section 5) on the product space $\mathcal{H}^n$, following the method of partial inverse of Spingarn (1983). Those are a few among numerous existing proximal splitting algorithms.

The increasing success of those methods during the past decades can be largely explained by their relative ease of implementation, their ability to handle large-scale problems, and their wide range of application. In particular, they allow to solve optimization problems when the monotone operators considered are subdifferentials of convex functionals. In this context, under suitable qualification conditions, solving (1.1) with our algorithm is equivalent to

$$\text{find } x \in \operatorname{argmin}\left\{\sum_{i=1}^{n} g_i + f\right\}, \tag{1.3}$$

by making only use of the *proximity operator* of each functional $g_i\colon \mathcal{H} \to ]-\infty, +\infty]$, supposed convex, proper and lower semicontinuous, and of the gradient of $f\colon \mathcal{H} \to \mathbb{R}$, supposed moreover *differentiable*, with a gradient which is *Lipschitz continuous* (equivalent to cocoercivity, by the theorem of Baillon and Haddad, 1977). Iteration (1.2) translates to this case by substituting $B$ with $\nabla f$ and each $J_{\frac{\gamma}{w_i}A_i}$ with $\mathrm{prox}_{\frac{\gamma}{w_i}g_i}\colon \mathcal{H} \to \mathcal{H}\colon x \mapsto \operatorname{argmin}_{y\in\mathcal{H}}\left\{\frac{1}{2}\|x-y\|^2 + \frac{\gamma}{w_i}g_i(y)\right\}$, the latter being the resolvent of the subdifferential of $\frac{\gamma}{w_i}g_i$ (Moreau, 1965).

In some cases however, first-order methods suffers from prohibitively slow rates of convergence. Consider in particular the simple gradient descent with constant step size, addressing (1.3) when each function $g_i$ is zero. Ignoring relaxation, iteration (1.2) reduces to $x \leftarrow x - \gamma\nabla f x$. If $\ell$ is (strictly positive) such that $\ell^{-1}\nabla f$ is *nonexpansive* (*i.e.* $\ell$ is a Lipschitz constant of $\nabla f$), then convergence towards a solution is ensured as long as $\ell\gamma < 2$.

However, when the gradient is significantly more sensitive along some directions than along others, *e.g.* for twice differentiable functionals, when the Hessian $\nabla^2 f$ is badly conditioned, a single scalar $\ell$ is a poor second-order information. A finer approach is to replace the (strictly positive) step size $\gamma$ by a (self-adjoint, strongly positive) linear operator $\Gamma$, adaptive to the direction: $x \leftarrow x - \Gamma\nabla f x$. If now $L$ is (self-adjoint, strongly positive) such that $L^{-1/2}\nabla f L^{-1/2}$ is nonexpansive, convergence condition becomes then $\|L^{1/2}\Gamma L^{1/2}\| < 2$.

Now, recall that if $M$ is any self-adjoint, strictly positive operator on $\mathcal{H}$, then $(x, y) \mapsto \langle Mx \mid y\rangle$ is an inner product on $\mathcal{H}$ and we denote $\mathcal{H}_M$ the Hilbert space $\mathcal{H}$ endowed with this inner product. Interestingly, the above term $\Gamma\nabla f$ is actually the gradient of $f$ in $\mathcal{H}_{\Gamma^{-1}}$. Similarly, $L^{-1/2}\nabla f L^{-1/2}$ being nonexpansive in $\mathcal{H}$ is equivalent to $f$ having a nonexpansive gradient in $\mathcal{H}_L$. Therefore, it is customary to refer to the above modification of the gradient descent as a *change of metric*.

Such method has been known for long (see for instance Davidon, 1959, 1991); in particular for twice differentiable $f$, setting at each iteration $\Gamma \stackrel{\text{set}}{=} \left(\nabla^2 f(x)\right)^{-1}$ yields *Newton method* for finding a zero of $\nabla f$. The quasi-Newton methods, where the inverse of the Hessian is iteratively estimated (see Broyden, 1967), are certainly the most celebrated variable metric methods for numerical optimization. Though they were originally developed for minimizing smooth functionals, variable metric methods can be applied to nonsmooth optimization techniques; the main purpose of this paper is to study an adapted metric version of our generalized forward-backward splitting algorithm. In the following, we give first an overview of preconditioning techniques for proximal algorithms. Then, we motivate our specific emphasis over applications to optimization structured on graphs. Finally, we precise our contributions and outline the content of this paper.

## 1.1 Preconditioning of Proximal Splitting Algorithms

Following ideas of Qian (1992), variable metric versions for the *proximal point algorithm* are provided in the context of optimization, for instance by Bonnans *et al.* (1995), Qi and Chen (1997) and Chen and Fukushima (1999). This allows to solve (1.3) in the restricted case $f \overset{\text{set}}{=} 0$ and only $n \overset{\text{set}}{=} 1$ functional (*i.e.* without splitting). They focus on convergence rates, usually involving some smoothness hypotheses. The proposed choice of the variable metric are derived from quasi-Newton methods, taking advantage of the fact that a proximal step on a functional is a gradient descent on its Moreau envelope; very little is said about the actual computation of the proximal step within the modified metric. Burke and Qian (1999) write the corresponding study of the proximal point algorithm in the general framework of monotone operators (solving (1.1) where $B \overset{\text{set}}{=} 0$, $n \overset{\text{set}}{=} 1$).

During the same period, Chen and Rockafellar (1997) give a variable metric version of the forward-backward splitting algorithm for monotone inclusion, addressing (1.1) in the restricted case $n \overset{\text{set}}{=} 1$, $A_1 + B$ strongly monotone and $\mathcal{H}$ finite-dimensional, achieving linear rate of convergence. The actual choice of the variable metric is not discussed, neither is the computation of the resolvent of $A_1$ (for the backward step) within the modified metric.

It is only a decade later that modified metric versions of proximal algorithm regained attention. Parente *et al.* (2008) explore respectively a variable metric inexact proximal point algorithm, and give a specific application to the resolution of a system of monotone differentiable equations. The same authors introduces then (Lotito *et al.*, 2009) a class of methods solving structured monotone inclusions (encompassing the particular case (1.1) with $n \overset{\text{set}}{=} 1$), retrieving some already existing methods as special case; no application is considered.

Pock and Chambolle (2011) present a preconditioning of their popular primal-dual proximal splitting algorithm (Chambolle and Pock, 2011) for minimizing functionals of the form $g_1 + g_2 \circ K$, where $K$ is any bounded linear operator, but which do not take advantage of the potential smoothness of a term. They do not consider a metric variable along iterations, hence the term *pre*conditioning. For computational purpose, they give a family of diagonal preconditioners which aim at reducing the condition number of the involved linear operators, and which are easy to compute as long as $K$ is known as a matrix.

Becker and Fadili (2012) focus on the computational point of view of the variable metric forward-backward splitting algorithm applied to convex optimization. They show how the proximal step of some functionals of interest can be efficiently computed when the metric can be decomposed as the sum of a diagonal operator and of an operator of rank 1. This allows them to consider the *symmetric rank* 1 quasi-Newton update (see Broyden, 1967), derived from the secant equation on $f$, *i.e.* on the smooth part of the splitting.

Chouzenoux *et al.* (2014) write a variable metric version of the inexact forward-backward of Attouch *et al.* (2013), able to address (1.3) in the special case $n \overset{\text{set}}{=} 1$ but with a possibly nonconvex smooth functional $f$. In numerical applications, they consider for $f$ a nonconvex data-fidelity term arising from many signal processing tasks, and show how to build metrics along iterations so as to obtain good majorizing, quadratic approximations. In order to compute the proximal step, which is in itself an optimization problem, they must however resort to subiterations.

Closer to our current interest, Giselsson and Boyd (2014a, 2014b) consider preconditioning of proximal methods for minimizing $f + g \circ K$ where $f$ and $g$ are convex, $K$ is a bounded linear operator and $f$ has additional properties. In the former paper, they assume $f$ strongly convex, which translates to a smooth functional in the splitting of the dual formulation, which is then solved with *fast forward-backward* splitting algorithm (see the review of Tseng, 2008). In the latter paper, they study the Douglas–Rachford splitting applied either to the primal or to the dual formulation, and derive an interesting linear convergence rate under the hypothesis of $f$ being both strongly convex and smooth. In both papers, they motivate the use of preconditioners which reduce the condition number of an operator which depends on $K$ and on the metrics capturing the properties of $f$. They also advise for diagonal preconditioners in order to ease the computations of proximity operators, and suggest various procedures to select them based on the above criterion.

In a more theoretical perspective, the first convergence proof of a variable metric forward-backward splitting algorithm is provided by Combettes and Vũ (2014), in the general framework of maximally monotone operators defined in Hilbert spaces. Vũ (2013) also gives the analogous study for the *forward-backward-forward* splitting algorithm of Tseng (2000), the latter being originally designed to tackle (1.1) in the restricted case $n \overset{\text{set}}{=} 1$, but with an operator $B$ which is merely Lipschitzian. Pesquet and Repetti (2014) write the preconditioning version of the *randomized block-coordinate* forward-backward splitting introduced recently by Combettes and Pesquet (2014). Unfortunately, none of these works exemplify how such variable metric can be used in practice. At last, Lorenz and Pock (2015) develop an *inertial forward-backward* splitting algorithm, allowing for a preconditioning and a variable step size along iterations. Among other practical considerations for convex optimization, they extend the diagonal preconditioners proposed by Pock and Chambolle (2011) to the case where the gradient of a smooth functional in the splitting are considered in explicit steps.

Let us recall that many problems can be reduced to the forward-backward setting when written in the primal-dual form over the appropriate augmented space. This also applies for problems whose scope go beyond (1.1) (or (1.3) in the context of convex optimization) by involving compositions with linear operators and parallel sums (inf-convolutions in convex optimization). Using this fact, each work mentioned in the previous paragraph develop variable metric versions of the primal-dual algorithms of Vũ (2013); Combettes and Pesquet (2012); Condat (2013). Although each of the latter methods can be used to solve (1.1) with full splitting, our generalized forward-backward, dealing only with the primal problem, is somewhat simpler and seems more adapted to splitting restricted to the form (1.1) (see Raguet *et al.*, 2013 and Raguet, 2014, § III.2.3, § IV.5). In this context, it is natural to write a preconditioned version of our algorithm.

## 1.2 Optimization on Graphs

As already pointed out, proximal splitting algorithms are mostly useful for solving problems where billions of coordinates are involved in an intricate way. A growing number of such problems are structured on *graphs*, either because it is physically inherent to the subject at hand (see for instance Shuman *et al.* (2013) for signal processing in general, and Richiardi *et al.* (2013) in the context of functional imaging), or because it captures abstract relationship between objects (see Peyré *et al.* (2011) in image processing).

Usually in this context, the splitting strategy aims at isolating some coordinates in the graph from interacting neighbors. A recent, notable example is provided by Couprie *et al.* (2013), who develop a dual optimization framework for tackling *total variation* based regularizations of inverse problems. The coordinates of interest are labels over vertices of the graph, and the proposed regularizations penalize finite differences (or increasing functions thereof) between neighboring vertices in the graph. The dimension of the dual space is thus the number of edges, which is usually much higher than the number of vertices; however, the dual formulation of some cases of interest can be easily split as the sum of simpler functionals.

In contrast in the primal-only setting, the splitting seems less obvious. Consider the simplest penalization inspired by total variation, which is the sum of the absolute values of the differences between all pairs of coordinates labelling adjacent vertices in the graph. The proximity operator of the distance between two scalars is easy to compute (see Proposition 4.1 (ii)), but in order to use proximal schemes like (1.2), each coordinate would have to appear within as many auxiliary variables as the degree (the number of neighbors) of its corresponding vertex. Thus, a naive application of a primal splitting would introduce, at least, as many auxiliary variables as the number of neighbors of the vertex of highest degree in the graph.

## 1.3 Tight Primal Splitting

To the best of our knowledge, in all primal-only splitting proposed to date, each auxiliary variable is a *full copy* of the original one. In the case described above, these auxiliary variables are supposed to have the dimension of the total number of vertices, regardless of the actual number of neighbors of the vertices of fewer degrees. This is not a concern when the graph is regular, but is prohibitively suboptimal in the case of highly irregular graphs. More generally, a naive primal splitting like (1.3) would not be *tight* if for some $i \in \{1, \dots, n\}$, the functional $g_i$ depends only on a strict subspace $\mathcal{H}_i \subset \mathcal{H}$, *i.e.* $g_i = g_i \circ P_{\mathcal{H}_i}$ where $P_{\mathcal{H}_i}$ is the orthogonal projector on $\mathcal{H}_i$.

Note that this situation can be cast in an even more general setting by considering any bounded linear operator $L_i \colon \mathcal{H} \to \mathcal{G}_i$, composed with $g_i \colon \mathcal{G}_i \to ]-\infty, +\infty]$. Then, the resulting problem can be efficiently tackled by the numerous primal-dual algorithms mentioned in § 1.1, where the auxiliary variables actually lie in the target Hilbert spaces $\mathcal{G}_i$. Now, it would be possible to extend our generalized forward-backward to this setting, provided that the operator $\sum_{i=1}^{n} w_i L_i{}^* L_i \colon \mathcal{H} \to \mathcal{H}$ is invertible, and at the expense of computing its inverse at each iteration. Although such extension does not require much effort and presents its own interest, it goes beyond the scope of this paper and is left for a future publication.

Instead, we limit the current study to the particular case of interest where for all $i \in \{1, \dots, n\}$, $L_i \overset{\text{set}}{=} P_{\mathcal{H}_i}$, and focus on reducing the auxiliary product space $\mathcal{H}^n$ to the restrained space $\times_{i=1}^{n} \mathcal{H}_i$. Let us already mention that this is illustrated in § 4.2.1 for the above example on a graph, where $\mathcal{H}$ is high-dimensional but most $\mathcal{H}_i$ lie in $\mathbb{R}^2$.

In the monotone inclusion setting (1.1), this case translates for each $i$ to $A_i = P_{\mathcal{H}_i} A_i P_{\mathcal{H}_i}$, that is to say, $A_i$ depends only on $\mathcal{H}_i$ but also has range included in $\mathcal{H}_i$. Moreover in the context of preconditioning, we replace the (strictly positive) weights $(w_i)_{1 \leq i \leq n}$ by suitable (self-adjoint, positive) linear operators $(W_i)_{1 \leq i \leq n}$. At the expense of adding a zero operator in the splitting (1.1), we can assume span $\left\{\cup_{i=1}^{n} \mathcal{H}_i\right\} = \mathcal{H}$, and the unit sum condition on the weights becomes $\sum_{i=1}^{n} W_i = \mathrm{Id}_{\mathcal{H}}$.

## 1.4 Contributions and Outline

We propose an extension of our generalized forward-backward splitting algorithm for solving monotone inclusion problems of the form (1.1), and in particular optimization problems of the form (1.3). This extension allows for preconditioning of the explicit forward step (through the operator $\Gamma$, generalizing the step size $\gamma$) as well as each of the implicit backward steps (through the operators $(W_i)_{1\leq i\leq n}$, generalizing the weights $(w_i)_{1\leq i\leq n}$). Moreover, it smartly takes into account situations where all the coordinates are not involved in the same numbers of operators in the splitting.

The necessary definitions and the algorithmic scheme is described in § 2, together with technical comparisons with other works. For the interested reader, the detailed convergence proof in an infinite dimensional setting is expounded in § 5. In § 3, we discuss practical preconditioning approaches for use of our generalized forward-backward for convex optimization. In particular, we provide some heuristics for dealing with popular regularization terms, such as sparsity-inducing norms and total variation.

Finally, motivated by problems encountered in geoinformatics, we exemplify these approaches in § 4 on popular regularizations of data represented over graphs. In particular, one of them is akin to total variation, as a convex relaxation of the geometric term in the Mumford-Shah functional. We show that the heuristics presented in § 3 give relevant preconditionings, and that the resulting minimization method compares favorably with alternatives from the literature.

## 2 Problem Formulation and Algorithm

We formulate the problem at hand together with our proposed algorithm, and then give the proof of convergence. For the sake of completeness, we first precise some definitions. Some properties of monotone operators are given in § 5.1; for general notions of convex analysis, see for instance Bauschke and Combettes (2011).

### 2.1 Preliminary Definitions

We consider a given real Hilbert space $\mathcal{H}$ with inner product $\langle\cdot\,|\,\cdot\rangle$ and associated norm $\|\cdot\|$. The identity operator on $\mathcal{H}$ is denoted $\mathrm{Id}_{\mathcal{H}}$ or simply Id when no confusion is possible. In the following, $A\colon \mathrm{dom}\,A \to 2^{\mathcal{H}}$ is a set-valued operator and $T\colon \mathrm{dom}\,T \to \mathcal{H}$ is a single-valued operator.

**Definition 2.1** (graph, inverse, domain, range, zeros and fixed-points)**.** The *graph* of $A$ is the set $\mathrm{gra}\,A \overset{\text{def}}{=} \{(x,y)\in\mathcal{H}^2 \mid y\in Ax\}$. The *inverse* of $A$, noted $A^{-1}$, is the operator whose graph is $\mathrm{gra}\,A^{-1} \overset{\text{def}}{=} \{(x,y)\in\mathcal{H}^2 \mid (y,x)\in\mathrm{gra}\,A\}$. The *domain* of $A$ is $\mathrm{dom}\,A \overset{\text{def}}{=} \{x\in\mathcal{H} \mid Ax\neq\emptyset\}$. The *range* of $A$ is $\mathrm{ran}\,A \overset{\text{def}}{=} \{y\in\mathcal{H} \mid \exists\, x\in\mathcal{H},\ y\in Ax\}$, and its set of *zeros* is $\mathrm{zer}\,A \overset{\text{def}}{=} \{x\in\mathcal{H} \mid 0\in Ax\} = A^{-1}0$. The set of *fixed-points* of $T$ is $\mathrm{fix}\,T \overset{\text{def}}{=} \{x\in\mathcal{H} \mid Tx = x\}$.

**Definition 2.2** (resolvent and reflection operators)**.** The *resolvent* of $A$ is the operator $J_A \overset{\text{def}}{=} (\mathrm{Id} + A)^{-1}$. The *reflection operator* associated to $J_A$ is the operator $R_A \overset{\text{def}}{=} 2J_A - \mathrm{Id}$.

**Definition 2.3** (maximally monotone operator)**.** $A$ is *monotone* in $\mathcal{H}$ if

$$\forall\,(x,y)\in(\mathrm{dom}\,A)^2,\quad \forall\,(u,v)\in Ax\times Ay,\quad \langle u-v \mid x-y\rangle \geq 0\,;$$

it is moreover *maximally monotone* if its graph is not strictly contained in the graph of any other monotone operator.

**Definition 2.4** (uniformly monotone operator)**.** $A$ is *uniformly monotone* in $\mathcal{H}$ with modulus $\phi\colon \mathbb{R}_+ \to [0,+\infty]$ if $\phi$ is a nondecreasing function vanishing only at 0 such that

$$\forall\,(x,y)\in(\mathrm{dom}\,A)^2,\quad \forall\,(u,v)\in Ax\times Ay,\quad \langle u-v \mid x-y\rangle \geq \phi(\|x-y\|)\,.$$

**Definition 2.5** (nonexpansive, $\alpha$-averaged and firmly nonexpansive operators)**.** $T$ is *nonexpansive* over $\mathcal{H}$ if

$$\forall\,(x,y)\in(\operatorname{dom}T)^2,\quad \|Tx-Ty\|\le\|x-y\|\,;$$

for $\alpha\in]0,1[$, it is moreover *$\alpha$-averaged* if there exists $R$ nonexpansive such that $T=\alpha R+(1-\alpha)\,\mathrm{Id}$. We denote $\mathcal{A}(\alpha,\mathcal{H})$ the set of all $\alpha$-averaged operators in $\mathcal{H}$. In particular, $\mathcal{A}(\frac{1}{2},\mathcal{H})$ is the set of *firmly nonexpansive* operators in $\mathcal{H}$.

**Definition 2.6** (cocoercive operator)**.** For $\beta\in]0,+\infty[$, $T$ is *$\beta$-cocoercive* in $\mathcal{H}$ if $\beta T\in\mathcal{A}\big(\frac{1}{2},\mathcal{H}\big)$.

Following for instance Minty (1962) and Zarantonello (1971), we know that an operator $A$ is maximally monotone over $\mathcal{H}$ if and only if its resolvent $J_A$ is firmly nonexpansive with full domain. According to Moreau (1965, Proposition 12.b), a particular case of interest of a maximally monotone operator is the subdifferential $\partial g$ of a proper, convex and lower semicontinuous $g\colon\mathcal{H}\to]-\infty,+\infty]$; Moreau (1965, Proposition 6.a) shows then that the resolvent $J_{\partial g}$ corresponds to the *proximity operator* $\operatorname{prox}_g\colon\mathcal{H}\to\mathcal{H}\colon x\mapsto\operatorname{argmin}_{y\in\mathcal{H}}\frac{1}{2}\|x-y\|^2+g(y)$.

Now if $C\subseteq\mathcal{H}$ is a closed convex subset of $\mathcal{H}$, we denote $P_C$ the *orthogonal projector* onto $C$; orthogonal projector are themselves particular cases of proximity operators. If moreover $\mathcal{S}\subseteq\mathcal{H}$ is a closed *subspace* of $\mathcal{H}$, recall that $P_{\mathcal{S}}$ is linear, idempotent and self-adjoint (see for instance Riesz and Sz.-Nagy, 1990, § 105); in that case, for any $x\in\mathcal{H}$, we use the convenient notation $x^{\mathcal{S}}\overset{\text{def}}{=}P_{\mathcal{S}}x$. We also denote the reflection $R_{\mathcal{S}}\overset{\text{def}}{=}2P_{\mathcal{S}}-\mathrm{Id}$.

Finally, we consider definitions mentioned in § 1.1. We say that a bounded, linear, self-adjoint operator $M\colon\mathcal{H}\to\mathcal{H}$ is *positive* if for all $x\in\mathcal{H}$, $\langle Mx\,|\,x\rangle\ge0$, and *strongly positive* if there exists a strictly positive real $m\in\mathbb{R}_{++}$ such that $M-m\,\mathrm{Id}$ is positive; we denote respectively by $\mathfrak{S}_+(\mathcal{H})$ and $\mathfrak{S}_{++}(\mathcal{H})$ the sets of such operators.

For $M\in\mathfrak{S}_+(\mathcal{H})$ (respectively $M\in\mathfrak{S}_{++}(\mathcal{H})$), $(x,y)\mapsto\langle x\,|\,y\rangle_M\overset{\text{def}}{=}\langle Mx\,|\,y\rangle$ and $x\mapsto\|x\|_M\overset{\text{def}}{=}\sqrt{\langle Mx\,|\,x\rangle}$ define a semi-inner product and a seminorm (respectively an inner product and a norm) over $\mathcal{H}$. For $M\in\mathfrak{S}_{++}(\mathcal{H})$, we denote by $\mathcal{H}_M$ the Hilbert space $\mathcal{H}$ endowed with this inner product. Since for all $x\in\mathcal{H}$, $m\|x\|^2\le\|x\|_M^2\le\|M\|\|x\|^2$, these norms are equivalent and induce the same topology over $\mathcal{H}$. In particular, if $\mathcal{S}\subseteq\mathcal{H}$ is closed in $\mathcal{H}$, then it is closed in $\mathcal{H}_M$, and we write $P_{\mathcal{S}}^M$ and $R_{\mathcal{S}}^M\overset{\text{def}}{=}2P_{\mathcal{S}}^M-\mathrm{Id}$ when considering an orthogonal projector over $\mathcal{S}$ in $\mathcal{H}_M$. More generally, if a functional $g$ is lower semicontinuous in $\mathcal{H}$, it is also lower semicontinuous in $\mathcal{H}_M$, and if it is also proper and convex, we write $\operatorname{prox}_g^M$ for its proximity operator in $\mathcal{H}_M$.

Finally, for $M\in\mathfrak{S}_+(\mathcal{H})$ (respectively $M\in\mathfrak{S}_{++}(\mathcal{H})$), we denote by $M^{1/2}$ (respectively $M^{-1}$) its *square root* (respectively *inverse*). We refer the reader to Riesz and Sz.-Nagy (1990, § 104) for the existence and properties of those operators.

## 2.2 Problem Statement and Assumptions

Let $n$ be a strictly positive integer. We consider problem (1.1) with the following assumptions.

(H1) $B\colon\mathcal{H}\to\mathcal{H}$ has full domain, and $L\in\mathfrak{S}_{++}(\mathcal{H})$ is such that

$$\forall\,(x,y)\in\mathcal{H}^2,\quad\langle Bx-By\,|\,x-y\rangle\ge\|Bx-By\|_{L^{-1}}^2\,. \tag{2.1}$$

(H2) For each $i\in\{1,\dots,n\}$, $A_i\colon\mathcal{H}\to2^{\mathcal{H}}$ is maximally monotone in $\mathcal{H}$, and $\mathcal{H}_i\subseteq\mathcal{H}$ is a closed subspace such that $A_i=P_{\mathcal{H}_i}A_iP_{\mathcal{H}_i}$.

(H3) $\operatorname{zer}\{\sum_{i=1}^{n} A_i + B\} \neq \emptyset$.

We now formulate the analogous assumptions in the convex optimization case (1.3).

(h1) $f : \mathcal{H} \to \mathbb{R}$ is convex and everywhere differentiable such that its gradient in $\mathcal{H}_L$ is nonexpansive, where $L$ is defined in (H1).

(h2) For each $i \in \{1, \ldots, n\}$, $g_i : \mathcal{H} \to ]-\infty, +\infty]$ is convex, proper and lower semicontinuous such that $g_i = g_i \circ P_{\mathcal{H}_i}$, where $\mathcal{H}_i$ is defined in (H2).

(h3) Domain qualification and feasibility conditions:

(i) for all $i \in \{1, \ldots, n\}$, $0 \in \operatorname{sri}\left\{\operatorname{dom} g_i - \cap_{j=1}^{i-1} \operatorname{dom} g_j\right\}$, and

(ii) $\operatorname{argmin}\{\sum_{i=1}^{n} g_i + f\} \neq \emptyset$.

Finally, we give the requirements on the preconditioners.

(P1) $\Gamma \in \mathfrak{S}_{++}(\mathcal{H})$ such that

(i) $\|L^{1/2} \Gamma L^{1/2}\| < 2$, where $L$ is defined in (H1), and

(ii) $\forall\, i \in \{1, \ldots, n\}$, $\Gamma(\mathcal{H}_i) \subseteq \mathcal{H}_i$, where $\mathcal{H}_i$ is defined in (H2).

(P2) For each $i \in \{1, \ldots, n\}$, $W_i \in \mathfrak{S}_{+}(\mathcal{H})$ is such that

(i) $\ker W_i = \mathcal{H}_i^{\perp}$,

(ii) $W_{i|\mathcal{H}_i} \in \mathfrak{S}_{++}(\mathcal{H}_i)$, and

(iii) $\Gamma^{-1} W_i = W_i \Gamma^{-1}$.

Moreover,

(iv) $\sum_{i=1}^{n} W_i = \mathrm{Id}_{\mathcal{H}}$.

**Remark 2.1.** Hypothesis (H2) is used to reduce the $i$-th auxiliary variable in Algorithm 1 from $\mathcal{H}$ to the possibly restricted space $\mathcal{H}_i$. Under this hypothesis, for all $i \in \{1, \ldots, n\}$, we consider $A_{i|\mathcal{H}_i}$, the restriction of $A_i$ to $\mathcal{H}_i$, *i.e.* the operator whose graph is $\operatorname{gra} A_{i|\mathcal{H}_i} \overset{\text{def}}{=} \operatorname{gra} A_i \cap {\mathcal{H}_i}^2$. For all $y \in \mathcal{H}_i$, thanks to Lemma 5.2 (i)⇔(iv) (see § 5.1) we can define $x \overset{\text{def}}{=} J_{A_i} y$. By definition, $y \in x + A_i x$, and since $A_i x \subseteq \mathcal{H}_i$, we deduce $x \in \mathcal{H}_i$, and thus $x = J_{A_{i|\mathcal{H}_i}} y$. Altogether, we get that $J_{A_{i|\mathcal{H}_i}} = (J_{A_i})_{|\mathcal{H}_i}$ is firmly nonexpansive in $\mathcal{H}_i$ with full domain, it is thus the resolvent of a maximally monotone operator in $\mathcal{H}_i$; but the latter is uniquely defined by $\left(J_{A_{i|\mathcal{H}_i}}\right)^{-1} - \mathrm{Id}_{\mathcal{H}_i} = A_{i|\mathcal{H}_i}$, proving maximal monotonicity of $A_{i|\mathcal{H}_i}$.

**Remark 2.2.** Denoting $\nabla f$ the gradient of $f$ in $\mathcal{H}$, it is easy to see that the gradient of $f$ in $\mathcal{H}_L$ is $L^{-1} \nabla f$. Thus,

$$
\begin{aligned}
\text{(h1)} &\Leftrightarrow \forall\,(x, y) \in \mathcal{H}^2,\ \|L^{-1}\nabla f x - L^{-1}\nabla f y\|_L \leq \|x - y\|_L\,, \\
&\Leftrightarrow \forall\,(x, y) \in \mathcal{H}^2,\ \|L^{-1/2}\nabla f x - L^{-1/2}\nabla f y\| \leq \|L^{1/2} x - L^{1/2} y\|\,, \\
&\Leftrightarrow \forall\,(u, v) \in \mathcal{H}^2,\ \|L^{-1/2}\nabla f L^{-1/2} u - L^{-1/2}\nabla f L^{-1/2} v\| \leq \|u - v\|\,,
\end{aligned}
$$

where we used the fact that $L^{1/2}$ is self-adjoint and invertible. We retrieve that (h1) is equivalent to nonexpansivity of $L^{-1/2}\nabla f L^{-1/2}$ in $\mathcal{H}$, as announced in the introduction. When $f$ is twice differentiable, a simple application of the chain rule and of the mean value inequality shows that a sufficient condition for (h1) is $\sup_{x \in \mathcal{H}} \{\|L^{-1/2}(\nabla^2 f(x)) L^{-1/2}\|\} \leq 1$.

**Remark 2.3.** Hypothesis (h2) is used to reduce the $i$-th auxiliary variable in Algorithm 1 in the convex optimization case from $\mathcal{H}$ to the possibly restricted space $\mathcal{H}_i$. Analogously to Remark 2.1, note that (h2) implies that $g_{i|\mathcal{H}_i}\colon \mathcal{H}_i \to \,]-\infty,+\infty]$ is convex, proper and lower semicontinuous in $\mathcal{H}_i$, and its proximity operator $\operatorname{prox}_{g_{i|\mathcal{H}_i}}$ is well defined.

**Remark 2.4.** A sufficient condition for (h3) (i) in finite dimension is $\cap_{i=1}^n \operatorname{ri}\{\operatorname{dom} g_i\} \neq \emptyset$. Also, condition (h3) (ii) is typically ensured by coercivity, or if any functional $g_i$ has bounded domain.

**Remark 2.5.** Referring again to Riesz and Sz.-Nagy (1990, § 104), under (P1) (ii) and (P2) (ii), (P2) (iii) imply that for all $i \in \{1,\dots,n\}$, $\Gamma^{-1}W_i \in \mathfrak{S}_{++}(\mathcal{H}_i)$.

**Algorithm 1** Preconditioned generalized forward-backward for solving monotone inclusion (1.1) under (H1)–(H3). For solving convex optimization (1.3) under (h1)–(h3), substitute $B$ with $\nabla f$ and for all $i \in \{1,\dots,n\}$, $J_{W_i^{-1}\Gamma A_i}$ with $\operatorname{prox}^{\Gamma^{-1}W_i}_{g_{i|\mathcal{H}_i}}$.

Require $(z_i)_{1\le i\le n} \in \times_{i=1}^n \mathcal{H}_i$ ; $\Gamma$ satisfying (P1) ; $(W_i)_{1\le i\le n}$ satisfying (P2) ;
$\forall\, k \in \mathbb{N}$, $\rho_k \in \,]0, 2 - \frac{1}{2}\|L^{1/2}\Gamma L^{1/2}\|[$ .

Initialization $x \leftarrow \sum_{i=1}^n W_i z_i$ ; $k \leftarrow 0$ .

**repeat**
  Main Iteration
  $p \leftarrow 2x - \Gamma B x$ ;
  **for** $i \in \{1,\dots,n\}$ **do**
    $z_i \leftarrow z_i + \rho_k\big(J_{W_i^{-1}\Gamma A_i}\big(p^{\mathcal{H}_i} - z_i\big) - x^{\mathcal{H}_i}\big)$ ;
  $x \leftarrow \sum_{i=1}^n W_i z_i$ ;
  $k \leftarrow k+1$ .
**until** *convergence*;
Return $x$.

## 2.3 Algorithmic Scheme

The algorithm for solving the monotone inclusion (1.1) or convex optimization (1.3) is given in Algorithm 1. The following theorem ensures the convergence and robustness to summable errors on the computations of each operator. For each iteration $k \in \mathbb{N}$, we define $b_k \in \mathcal{H}$ the error when computing $B$ and for each $i \in \{1,\dots,n\}$, $a_{i,k} \in \mathcal{H}_i$ the error when computing $J_{W_i^{-1}\Gamma A_i}$.

**Theorem 2.1.** *Set $(z_{i,0})_{1\le i\le n} \in \times_{i=1}^n \mathcal{H}_i$ and define $(x_k)_{k\in\mathbb{N}}$ the sequence in $\mathcal{H}$ together with $\big((z_{i,k})_{1\le i\le n}\big)_{k\in\mathbb{N}}$ the sequence in $\times_{i=1}^n \mathcal{H}_i$ such that for all $k \in \mathbb{N}$, $x_k = \sum_{i=1}^n W_i z_{i,k}$ and for all $i \in \{1,\dots,n\}$,*

$$z_{i,k+1} = z_{i,k} + \rho_k\Big(J_{W_i^{-1}\Gamma A_i}\big(2x_k{}^{\mathcal{H}_i} - z_{i,k} - \Gamma(Bx_k + b_k)^{\mathcal{H}_i}\big) + a_{i,k} - x_k{}^{\mathcal{H}_i}\Big)\,, \tag{2.2}$$

*where $b_k \in \mathcal{H}$, $a_{i,k} \in \mathcal{H}_i$, and $\rho_k \in \,]0, 2 - \frac{1}{2}\|L^{1/2}\Gamma L^{1/2}\|[$.*
*Under assumptions (H1)–(H3) and (P1)–(P2), if*

(i) $\sum_{k\in\mathbb{N}} \rho_k\big(2 - \frac{1}{2}\|L^{1/2}\Gamma L^{1/2}\| - \rho_k\big) = +\infty$, *and*

(ii) $\sum_{k\in\mathbb{N}} \rho_k\|b_k\| < +\infty$, *and* $\forall\, i \in \{1,\dots,n\}$, $\sum_{k\in\mathbb{N}} \rho_k\|a_{i,k}\| < +\infty$,

*then the sequence* $(x_k)_{k\in\mathbb{N}}$ *defined by* (2.2) *converges weakly towards a solution of* (1.1). *If moreover*

(iii) $0 < \inf_{k\in\mathbb{N}} \{\rho_k\} \le \sup_{k\in\mathbb{N}} \{\rho_k\} \le 2 - \frac{1}{2}\|L^{1/2}\Gamma L^{1/2}\|$, *and*

(iv) $B$ *is uniformly monotone,*

*then the convergence is strong.*

The following corollary specializes Theorem 2.1 to the case of convex optimization problems of the form (1.3).

**Corollary 2.1.** *Let the sequence* $(x_k)_{k\in\mathbb{N}}$ *be defined by substituting in* (2.2), $B$ *with* $\nabla f$ *and for all* $i \in \{1,\dots,n\}$, $J_{W_i^{-1}\Gamma A_i}$ *with* $\mathrm{prox}^{\Gamma^{-1}W_i}_{g_i|_{\mathcal{H}_i}}$. *If, in addition to* (h1)–(h3) *and* (P1)–(P2), *assumptions* (i)–(ii) *of* Theorem 2.1 *are satisfied, then the sequence* $(x_k)_{k\in\mathbb{N}}$ *converges weakly towards a minimizer of* (1.3). *If moreover assumption* (iii) *of* Theorem 2.1 *is satisfied, and* $f$ *is uniformly convex, then* $(x_k)_{k\in\mathbb{N}}$ *converges strongly to the unique minimizer of* (1.3).

**Remark 2.6** (Scalar metric). Setting $\Gamma \overset{\text{set}}{=} \gamma\,\mathrm{Id}_{\mathcal{H}}$, $\gamma \in \mathbb{R}_{++}$, and for all $i \in \{1,\dots,n\}$, $\mathcal{H}_i \overset{\text{set}}{=} \mathcal{H}$ and $W_i \overset{\text{set}}{=} w_i\mathrm{Id}_{\mathcal{H}}$, $w_i \in \mathbb{R}_{++}$, such that $\sum_{i=1}^n w_i = 1$, one retrieves our generalized forward-backward splitting algorithm. Note that condition (H1) (respectively condition (h1)) with $L \overset{\text{set}}{=} \ell\,\mathrm{Id}_{\mathcal{H}}$, $\ell \in \mathbb{R}_{++}$, is equivalent to $B$ being cocoercive with modulus $\ell^{-1}$ (respectively to $\nabla f$ being Lipschitz continuous with constant $\ell$). In those cases, (P1) (i) reduces to the classical convergence condition on the step size $\gamma\ell < 2$, and condition (i) on the relaxation parameters becomes $\sum_{k\in\mathbb{N}} \rho_k\big(2 - \frac{1}{2}\gamma\ell - \rho_k\big) = +\infty$ as given by us in its full extent elsewhere (Raguet, 2014, Theorem III.2.1 (ii)). Finally, condition (iii) in Theorem 2.1 allows strong convergence from uniform monotonicity of $B$ with overrelaxation, *i.e.* $1 < \sup_{k\in\mathbb{N}} \{\rho_k\} < 2 - \frac{\gamma\ell}{2}$. Shortly before publishing our work, we became aware of the recent work of Combettes and Yamada (2015), proposing independently and concurrently the same bound for the simple forward-backward splitting.

**Remark 2.7** (Variable metric forward-backward splitting). In the special case $n \overset{\text{set}}{=} 1$, our algorithm reduces to a preconditioned version of the forward-backward splitting. In this case, it is possible to use variable metrics, *i.e.* varying preconditioning matrix $\Gamma$ along iterations, as studied by Combettes and Vũ (2014) and Chouzenoux *et al.* (2014), and to a certain extent by Lorenz and Pock (2015) who only consider varying step sizes. Also in the context of convex optimization, accelerated versions of the forward-backward can be easily preconditioned, as considered by Giselsson and Boyd (2014a). Finally, note that the definition of $L$, the cocoercivity metric of $B$, given by (2.1), and the resulting condition (P1) (i) on the step metric $\Gamma$, are somewhat similar to those given by Lorenz and Pock (2015) or Pesquet and Repetti (2014). In contrast, Combettes and Vũ (2014) restrict their analysis to a scalar metric $L$, yielding possible step metrics which are of little interest for practical preconditioning purpose.

**Remark 2.8** (Preconditioned Douglas–Rachford splitting). In the special case $B \overset{\text{set}}{=} 0$, $n \overset{\text{set}}{=} 2$ and for all $i \in \{1,\dots,n\}$, $\mathcal{H}_i \overset{\text{set}}{=} \mathcal{H}$, our algorithm reduces to a preconditioned version of the relaxed Douglas–Rachford splitting. This setting is studied in Giselsson and Boyd (2014b), in both the primal and dual formulations. For $n \le 2$, our method constitutes a preconditioning of the Douglas–Rachford splitting applied to the method of the partial inverses on the product space $\mathcal{H}^n$ (see § 1). In addition, it allows for more adapted splitting along coordinates with arbitrary $\mathcal{H}_i \subset \mathcal{H}$ as explained in § 1.3, and for individual preconditioning of each function in the splitting, as further detailed in § 3.2.

# 3 Preconditioning Approaches for Primal Convex Optimization

We detailed in the previous section the necessary conditions over the involved operators to ensure convergence of our preconditioned generalized forward-backward. In this section, we describe some approaches for choosing the preconditioners in practice, in the context of convex optimization.

As mentioned in the introduction § 1.2, our specific interest lies on large-scale problems, involving sums of functionals that are fairly simple when taken individually. For reasons which will become clear in the following, we restrict ourselves to diagonal preconditioners.

With this focus in mind, we first exemplify the usefulness of preconditioning on specific but popular instances of functionals, before describing a general method for computing the preconditioners for our generalized forward-backward splitting algorithm.

## 3.1 First-order Primal Proximal Preconditioning

We briefly present here two of the main motivations for preconditioning first order proximal methods. First, we show that preconditioning can drastically simplify the computation of proximity operators. Then, we describe the use of second-order information for accelerating the convergence, following classic changes of metric inspired by quasi-Newton methods (see introduction).

### 3.1.1 Simplifying some Proximity Operators

The most straightforward examples are convex penalizations involving the norm of the considered Hilbert space; let us mention some of them. In the terms and notations of the previous sections, let $\mathcal{S} \subseteq \mathcal{H}$ be a particular closed subspace of interest, let $\lambda \in \mathbb{R}_{++}$, and define the *group seminorm* $g\colon x \mapsto \lambda \|x^{\mathcal{S}}\|$, and the *group constraint*

$$h\colon x \mapsto \begin{cases} +\infty & \text{if } \|x^{\mathcal{S}}\| > \lambda\,, \\ 0 & \text{if } \|x^{\mathcal{S}}\| \leq \lambda\,. \end{cases}$$

Note that when $\mathcal{S} \overset{\text{set}}{=} \mathcal{H}_i$ is spanned by a subset of a given orthonormal basis of $\mathcal{H}$, $g$ is an example of *group norm* over this basis; in this case, it is often used as a regularizing functional to enforce *group sparsity* (see for instance Jenatton *et al.*, 2011) in the solution of (1.3). Another case of interest is when $\mathcal{S} = {\mathcal{D}_i}^{\perp} \cap \mathcal{H}_i$, where $\mathcal{D}_i$ is the *first bisector* in a given orthonormal basis of $\mathcal{H}_i$. The functional $g$ is then the *deviation seminorm*, which can be used to discretize the total variation seminorm with better isotropy and more convenient computational properties than classical discretizations (see Raguet, 2014, § IV.1.3, § IV.2).

In any case, the proximity operators of these functionals only involve the computations of the projector over $\mathcal{S}$ and of the norm, as shown in the following proposition.

**Proposition 3.1.** *With the above definitions, for all $x \in \mathcal{H}$,*

$$\begin{cases} \textit{if } \|x^{\mathcal{S}}\| > \lambda, \textit{ then } \operatorname{prox}_g(x) = x^{\mathcal{S}^\perp} + \left(1 - \frac{\lambda}{\|x^{\mathcal{S}}\|}\right) x^{\mathcal{S}} \textit{ and } \operatorname{prox}_h(x) = x^{\mathcal{S}^\perp} + \frac{\lambda}{\|x^{\mathcal{S}}\|} x^{\mathcal{S}}\,; \\ \textit{if } \|x^{\mathcal{S}}\| \leq \lambda, \textit{ then } \operatorname{prox}_g(x) = x^{\mathcal{S}^\perp} \qquad\qquad \textit{and } \operatorname{prox}_h(x) = x\,. \end{cases}$$

*Proof.* For $\mathcal{S} \overset{\text{set}}{=} \mathcal{H}$, $\operatorname{prox}_g$ and $\operatorname{prox}_h$ reduces respectively to the *group soft-thresholding* (Combettes and Wajs, 2005, Example 2.16) and to the orthogonal projector over the closed ball centered in 0 and of radius $\lambda$, corresponding to the above expressions. Using this fact, together with the fact that $P_S$ has the *tight frame property* with squared norm 1, *i.e.* $P_S{}^* P_S = P_S$, an earlier result on compositions of functionals with such operators (see Raguet, 2014, Proposition IV.3.7) applies and lead to the given expressions. ■

Now in many practical circumstances, it makes sense to use weighted versions of group seminorms or group constraints, for instance to account for different significations or importances along the coordinates. This is easily enforced with a change of metric; letting $M \in \mathfrak{S}_{++}(\mathcal{H})$, we modify $g$ and $h$ above by considering the orthogonal projector over $\mathcal{S}$ in the metric induced by $M$, *i.e.* $g^M : x \mapsto \lambda \|P_S^M x\|_M$ and

$$h^M : x \mapsto \begin{cases} +\infty & \text{if } \|P_S^M x\|_M > \lambda\,, \\ 0 & \text{if } \|P_S^M x\|_M \le \lambda\,. \end{cases}$$

Usually, such $M$ is a diagonal operator, and the computations of the projector over $\mathcal{S}$ and of the norm in $\mathcal{H}_M$ are as easy as their scalar metric counterparts. Since no assumption is made over the scalar product in $\mathcal{H}$ in Proposition 3.1, we directly obtain the following corollary.

**Corollary 3.1.** *With the above definitions, for all $x \in \mathcal{H}$,*

$$\begin{cases} \textit{if } \|P_S^M x\|_M > \lambda, \textit{ then } \operatorname{prox}^M_{g^M}(x) = P^M_{S^\perp} x + \Big(1 - \frac{\lambda}{\|P_S^M x\|_M}\Big) P_S^M x \\ \qquad\qquad \textit{and } \operatorname{prox}^M_{h^M}(x) = P^M_{S^\perp} x + \frac{\lambda}{\|P_S^M x\|_M} P_S^M x\,; \\ \textit{if } \|P_S^M x\|_M \le \lambda, \textit{ then } \operatorname{prox}^M_{g^M}(x) = P^M_{S^\perp} x \textit{ and } \operatorname{prox}^M_{h^M}(x) = x\,. \end{cases}$$

**Remark 3.1.** Note that in our setting, we are actually interested in $M_i \in \mathfrak{S}_{++}(\mathcal{H}_i)$ and want to compute the proximity operator of $g^{M_i}{}_{|\mathcal{H}_i}$ within $\mathcal{H}_{i\,M_i}$. The substitution is straightforward.

In contrast the minimization problems defining $\operatorname{prox}_{g^M}$ or $\operatorname{prox}_{h^M}$ within the original metric can be much more complicated, precisely because it involves two different norms. Think for instance of the orthogonal projection over an ellipsoid, for which no general closed-form expression is available. Thus, the change of metric allows to take efficiently into account functionals in (1.3) such as weighted norms, weighted total variation and ellipsoid constraints.

### 3.1.2 Accelerating First-order Methods

As pointed out in the introduction, variable metrics for accelerating first-order methods based on gradient descent have been used for long and are well studied. In contrast, variable metric acceleration of methods involving proximity operators is not quite yet established. Recall from § 1.1 that attempts in that directions usually rely on some smoothness hypothesis; either for the proximal point algorithm (Bonnans *et al.*, 1995; Qi and Chen, 1997; Chen and Fukushima, 1999) or for the forward-backward splitting (Becker and Fadili, 2012; Chouzenoux *et al.*, 2014). In particular, note that the latter two works only use second-order information of $f$, the smooth part in the splitting. At last, recall the linear convergence rate of the Douglas–Rachford splitting algorithm obtained by Giselsson and Boyd (2014b) under the hypothesis that one functional in the splitting is both smooth and strongly convex.

Roughly speaking, functionals which are both smooth and strongly convex are close to *quadratic functionals*, $q : x \mapsto \langle Mx + \mu \,|\, x\rangle + q(0)$, with $M \in \mathfrak{S}_{++}(\mathcal{H})$ and $\mu \in \mathcal{H}$. In that case, the result of Giselsson and Boyd (2014b, Proposition 6) states that the rate is determined by the *condition number* of $M$, *i.e.* the ratio between its greatest and its smallest eigenvalues. Accordingly, this suggests to set the algorithm in the metric of $\mathcal{H}_M$, which obviously reminds of Newton method, since in that case the Hessian $\nabla^2 q$ is constant and equal to $M$.

Thus, our approach is quite natural: find a *quadratic approximation* of each functional involved in the splitting (1.3), and use the linear operator of the resulting quadratic form for preconditioning. Note that our approach is only heuristic, in order to speed up the convergence of the algorithm, and should not be confused with majorization-minimization techniques. In particular, the approximations need not be majorizers of the corresponding functionals.

We do not discuss how to obtain convenient quadratic approximations, neither how to ensure that the preconditioner is strongly positive and, as needed for most cases of interest, diagonal. Concerning the two last points, an interesting discussion is given by Giselsson and Boyd (2014a, § 6). We provide specific examples in § 4.2 and emphasize that for such large-scale, simple functionals, the diagonal approximation is obtained by simply dropping the off-diagonal terms (in the terms of the above cited discussion, this is coined *infinite norm equilibration*).

Finally, nothing has been said on how to combine several such metrics for dealing with an arbitrary splitting of the form (1.3); this is the subject of the next section.

## 3.2 Preconditioning our Generalized Forward-Backward Splitting

Suppose that we are given the preconditioners $M$ and $(M_i)_{1 \le i \le n}$, suggested respectively by the functionals $f$ and $(g_i)_{1 \le i \le n}$, either for computational (see § 3.1.1) or acceleration (see § 3.1.2) purposes. By construction, $M \in \mathfrak{S}_{++}(\mathcal{H})$ and for all $i \in \{1, \ldots, n\}$, $M_i \in \mathfrak{S}_{+}(\mathcal{H})$ such that $M_{i|\mathcal{H}_i} \in \mathfrak{S}_{++}(\mathcal{H}_i)$ and $\ker M_i = \mathcal{H}_i{}^{\perp}$.

Suppose moreover that there exists an orthogonal basis of $\mathcal{H}$ indexed by a set $J$, which contains for each $i \in \{1, \ldots, n\}$, an orthogonal basis of $\mathcal{H}_i$ indexed by $J_i \subseteq J$, and such that each of the above operators is diagonal in this basis. Although not theoretically necessary, this greatly simplifies the manipulations presented below. For convenience, denote $M \overset{\text{def}}{=} \operatorname{diag}(m_j)_{j \in J}$ and for all $i \in \{1, \ldots, n\}$, $M_i \overset{\text{def}}{=} \operatorname{diag}(m_{ij})_{j \in J}$; note in particular that $J_i = \{j \in J \,|\, m_{ij} > 0\}$.

### 3.2.1 Satisfying the Requirements

Looking at Algorithm 1 with § 3.1 in mind, one would like to precondition the algorithm with help of $\Gamma$ and $(W_i)_{1 \le i \le n}$ such that $\Gamma^{-1} \approx M$ and for all $i \in \{1, \ldots, n\}$, $\Gamma^{-1} W_i \approx M_i$. Unfortunately, in general we can't define them directly in this way, because requirements (P1)–(P2) won't be satisfied. We present here practical approaches for complying to these requirements.

**Preliminary: select the relaxation parameters $(\rho_k)_{k \in \mathbb{N}}$.** For ease of demonstration, we suppose that the relaxation parameters are set prior to the conditioners. Looking at Theorem 2.1, we have necessarily for all $k \in \mathbb{N}$, $\rho_k \in \left]0, 2\right[$. Many papers advocate the use of overrelaxation, so in the absence of other requirements, we set typically for all $k \in \mathbb{N}$, $\rho_k \overset{\text{set}}{=} 1.5$. Strictly speaking, it is even possible to set $(\rho_k)_{k \in \mathbb{N}}$ converging to a value arbitrarily close to 2, but it should be noted from (3.2) below that high relaxations impede the value of $\Gamma$.

Note that some specific implementations (see for intance Raguet, 2014, § IV.4.1) are more efficient without relaxation (for all $k \in \mathbb{N}$, $\rho_k \overset{\text{set}}{=} 1$). In some extreme setting where the operators are computed up to nonsummable errors, it might also be useful to consider $(\rho_k)_{k \in \mathbb{N}}$ converging to 0 in such a way that conditions (i) and (ii) of Theorem 2.1 are both satisfied; very slow convergence should be expected however.

**First step: creating the global preconditioner $\Gamma$.** A first approximation of the global preconditioner $\tilde{\Gamma} \overset{\text{def}}{=} \operatorname{diag}(\tilde{\gamma}_j)_{j\in J}$ can be immediately obtained either by preconditioning according to $f$ alone, or, for reasons that can be understood below, according to the whole functional. That is,

$$\tilde{\Gamma} \overset{\text{set}}{=} M^{-1}, \quad \textit{i.e.} \quad \forall\, j \in J,\ \tilde{\gamma}_j \overset{\text{set}}{=} 1/m_j\,, \tag{3.1a}$$

or

$$\tilde{\Gamma} \overset{\text{set}}{=} \left(M + \textstyle\sum_{i=1}^{n} M_i\right)^{-1}, \quad \textit{i.e.} \quad \forall\, j \in J,\ \tilde{\gamma}_j \overset{\text{set}}{=} 1\Big/\left(m_j + \textstyle\sum_{i=1}^{n} m_{ij}\right)\ . \tag{3.1b}$$

Note now that any diagonal operator satisfies condition (P1) (ii). In order to fulfill requirement (P1) (i), suppose that assumption (h1) is also given with a diagonal operator $L \overset{\text{def}}{=} \operatorname{diag}(\ell_j)_{j\in J}$, and simply set $\Gamma \overset{\text{def}}{=} \operatorname{diag}(\gamma_j)_{j\in J}$, with for all $j \in J$,

$$\gamma_j \overset{\text{set}}{=} \min\left\{\delta\frac{4-2\bar{\rho}}{\ell_j}, \tilde{\gamma}_j\right\}, \tag{3.2}$$

where $\bar{\rho} \overset{\text{def}}{=} \sup_{k\in\mathbb{N}}\{\rho_k\}$ and $\delta$ is some value less than 1, typically $\delta \overset{\text{def}}{=} 0.99$.

**Second step: creating the preconditioners $(W_i)_{1\le i\le n}$.** We have now set $\Gamma$, so that for all $i \in \{1,\dots,n\}$, a first approximation of the preconditioner $\tilde{W}_i \overset{\text{def}}{=} \operatorname{diag}(\tilde{w}_{ij})$ is given directly by $\tilde{W}_i \overset{\text{set}}{=} \Gamma M_i$, *i.e.* for all $j \in J$, $\tilde{w}_{ij} \overset{\text{set}}{=} \gamma_j m_{ij}$. By construction, each of these preconditioners satisfy conditions (P2) (i)–(iii). However, condition (P2) (iv) is not satisfied in general, that is to say $\tilde{S} \overset{\text{def}}{=} \sum_{i=1}^{n} \tilde{W}_i \neq \mathrm{Id}$. Observe yet that $\tilde{S} = \operatorname{diag}(\tilde{s}_j)_{j\in J}$, where for all $j \in J$, $\tilde{s}_j \overset{\text{def}}{=} \sum_{i=1}^{n} \tilde{w}_{ij}$. In order to fulfill the last condition, we propose two different ways of modifying each $\tilde{W}_i$, depending on the importance of keeping the shape of the metrics they induce.

Indeed, if it is not crucial to keep the exact metrics, as is typically the case for acceleration purpose (see again § 3.1.2), then the simplest way is to scale each coordinate separately, according to the relative amplitude of each preconditioner. That is, for all $i \in \{1,\dots,n\}$, set $W_i \overset{\text{def}}{=} \operatorname{diag}(w_{ij})_{j\in J}$ as

$$W_i \overset{\text{set}}{=} \tilde{S}^{-1}\tilde{W}_i, \quad \textit{i.e.} \quad \forall\, j \in J,\ w_{ij} \overset{\text{set}}{=} \tilde{w}_{ij}/\tilde{s}_j\ . \tag{3.3a}$$

Now, suppose in contrast that the shape of the metrics is indeed crucial, for instance for computational purpose as explained in § 3.1.1. Then, the only modifications allowed to the above preconditioners are multiplications by (positive) scalars. In particular, for all $i \in \{1,\dots,n\}$, define $s_i \overset{\text{def}}{=} \sup\{\tilde{s}_j \,|\, j \in J_i\}$, and set

$$W_i \overset{\text{set}}{=} \tfrac{1}{s_i}\tilde{W}_i\,, \quad \textit{i.e.} \quad \forall\, j \in J,\ w_{ij} \overset{\text{set}}{=} \tilde{w}_{ij}/s_i\ . \tag{3.3b}$$

Then, defining $W_{n+1} \overset{\text{def}}{=} \mathrm{Id} - \sum_{i=1}^{n} W_i = \operatorname{diag}\left(1 - \sum_{i=1}^{n} w_{ij}\right)_{j\in J}$, $\mathcal{H}_{n+1} \overset{\text{def}}{=} (\ker W_{n+1})^{\perp}$ and $g_{n+1} \overset{\text{def}}{=} 0$, it is easy to check that, substituting $n$ with $n+1$, all assumptions are fulfilled to apply Theorem 2.1, and that the sequence generated by the algorithm converges towards a solution of (1.3).

**Remark 3.2.** The proposed way of defining the scalars $(s_i)_{1 \le i \le n}$ is optimal in the sense that if any of them is decreased, then $W_{n+1}$ is not positive anymore. Let us precise also that the above technique applies of course to the monotone inclusion setting, substituting $g_{n+1}$ with $A_{n+1} \overset{\text{def}}{=} 0$ if necessary. In any case, the proximity operator of $g_{n+1}$, or the resolvent of $A_{n+1}$, reduce to the identity operator. In general, it is necessary to store in memory the auxiliary variable $z_{n+1}$, which has the dimension of the full problem; note however that using no relaxation (for all $k \in \mathbb{N}$, $\rho_k \overset{\text{set}}{=} 1$), this is not needed anymore since its update in Algorithm 1 reduces to $z_{n+1} \leftarrow p - x$. Finally, note that setting $\Gamma$ according to (3.1b) rather than (3.1a) is justified by the fact that $\Gamma$ intervenes in all final metrics. This was actually more efficient for all numerical experiments of § 4 (in which we only show results with method (3.1b)).

### 3.2.2 Stationary Metric and Reconditioning

We advocate in § 3.1.2 to select the preconditioning metrics according to quadratic approximations of the involved functionals. A major issue is that those approximations are only *local*, in the sense that the error is only controlled in a bounded neighborhood. Classically in such case, the approximations must be updated along with the iterate. Accordingly, the preconditioning operators must then be updated (following § 3.2.1) along the run of the algorithm; we coin such a process *reconditioning* step. Unfortunately, our preconditioned generalized forward-backward splitting does not allow for variable metric as soon as $n > 1$, so that some special care must be taken.

For better understanding, let us emphasize the theoretical problem by briefly referring to § 5.3. Given the preconditioners $\Gamma$ and $(W_i)_{1 \le i \le n}$, Proposition 5.2 ensures the existence of auxiliary points satisfying the fixed-point equation (5.2). The auxiliary variables $(z_i)_{1 \le i \le n}$ introduced in Algorithm 1 are indeed supposed to converge to such points. Now, if we modify the preconditioners, we also modify these auxiliary points. In particular, at the very next iteration following a reconditioning, the former $(z_i)_{1 \le i \le n}$ are further away from the new convergence points, so that one might expect the iterate $x$ to also draw away from the solution set.

In order to avoid such drawback, we thus suggest to also update the auxiliary variables when reconditioning, accordingly with the new preconditioners. A fairly simple approach is to suppose that convergence is almost achieved and that (5.2) is almost satisfied. Having a look at the proof of Proposition 5.2, this means that for all $i \in \{1, \ldots, n\}$, there exists $y_i \in A_i x$ such that $W_i(x - \Gamma B x - z_i) \approx \Gamma y_i$. Given the *updated* preconditioners $\hat{\Gamma}$ and $(\hat{W}_i)_{1 \le i \le n}$, we keep the same iterate $x$ and update the auxiliary variables $(\hat{z}_i)_{1 \le i \le n}$ according to, for all $i \in \{1, \ldots, n\}$,

$$\hat{z}_i = (x - \hat{\Gamma} B x)^{\mathcal{H}_i} - \hat{W}_i^{-1} \hat{\Gamma} y_i \;, \quad \text{where} \quad y_i = \Gamma^{-1} W_i (x - \Gamma B x - z_i) \;.$$

Note that this can be easily implemented within Algorithm 1 without significant computational or memory overhead.

In the spirit of Liang *et al.* (2014, Theorem 6.4), we could consider a variable metric algorithm as a stationary metric one with additional errors, and check under which conditions over the variable metrics the resulting errors are summable. We leave however such a study for a future work, and opt instead for a finite number of reconditionings along the run of the algorithm.

When the stopping criterion of the algorithm is specified as a maximum number of iterations, it is reasonable to perform reconditionings after given fractions of this number have been reached. Another more adaptive way is to monitor the relative evolution of the iterate, and perform reconditionings when this evolution goes below a certain threshold.

We do not extensively investigate how to optimally schedule the reconditionings. Let us only underline that, early in the optimization process, as the iterate is still far from the solution set, neither the quadratic approximations nor the condition (5.2) might be relevant. On the other hand, too late a reconditioning is futile, as its purpose is precisely to speed up the convergence. This is illustrated in our numerical experiments, § 4.3.1.

# 4 Application to Graph-Structured Optimization Problems

We now present a numerical application of our algorithm on a high-dimensional problem, structured on an irregular graph. Solving it implies performing a sequence of badly conditioned, nondifferentiable optimization problems, providing a good illustration of several notions developed in the previous sections. We first introduce the motivations and the terms of the problem. Then, we describe in details how we apply our preconditioned generalized forward-backward splitting. Finally, we show how the convergence speed of our approach compare with concurrent preconditioned proximal algorithms over three different instances of the problem, and provide visual illustrations of the results.

## 4.1 Aggregating Spatial Statistics

The great amount of georeferenced socio-economic data available exceeds our current ability to analyse them. We consider in another work (Landrieu *et al.*, 2015) the problem of aggregating spatial statistics to obtain simple yet accurate representations in map form, facilitating analysis and providing with valuable decision aids. The spatial data consists of an intensive value (*e.g.* the result of an election in percentages), related to an extensive quantity (*e.g.* the number of voters), defined over subregions of a geographical space (*e.g.* electoral constituencies). See some instances in Table 1.

In order to capture the spatial structure of the data, we consider a graph $\mathsf{G} \overset{\text{def}}{=} (\mathsf{V}, \mathsf{E})$, where the vertices $\mathsf{V}$ represent the subregions and the edges $\mathsf{E} \subseteq \mathsf{V} \times \mathsf{V}$ represent spatial adjacencies. Further spatial information is available, encoded by strictly positive vectors $\lambda \overset{\text{def}}{=} (\lambda_\mathsf{v})_{\mathsf{v}\in\mathsf{V}} \in \mathbb{R}_{++}^{|\mathsf{V}|}$ and $\mu \overset{\text{def}}{=} (\mu_{\mathsf{uv}})_{(\mathsf{u},\mathsf{v})\in\mathsf{E}} \in \mathbb{R}_{++}^{|\mathsf{E}|}$, weighting the vertices (respectively the edges) by their corresponding surface (respectively border length). Observed spatial data constitute vectors $y \overset{\text{def}}{=} (y_\mathsf{v})_{\mathsf{v}\in\mathsf{V}} \in \mathbb{R}^{|\mathsf{V}|}$ and $\nu \overset{\text{def}}{=} (\nu_\mathsf{v})_{\mathsf{v}\in\mathsf{V}} \in \mathbb{R}^{|\mathsf{V}|}$ labelling the vertices respectively by their intensive value and extensive quantity. In addition, our setting takes into account vertices $\mathsf{v}$ corresponding to regions for which either the value of $y_\mathsf{v}$ is missing, or the quantity $\nu$ is zero or undefined. In the following, we set by convention for such vertex $y_\mathsf{v} \overset{\text{set}}{=} 0$ or $\nu_\mathsf{v} \overset{\text{set}}{=} 0$ if necessary.

### 4.1.1 Problem Formulation

The simplification of the observed data $y$ is modeled as an element $x^{(0)} \in \operatorname{argmin}\{F_0\}$, where for all $x \overset{\text{def}}{=} (x_\mathsf{v})_{\mathsf{v}\in\mathsf{V}} \in \mathbb{R}^{|\mathsf{V}|}$,

$$F_0(x) \overset{\text{def}}{=} \frac{1}{2}\sum_{\mathsf{v}\in\mathsf{V}} \lambda_\mathsf{v}^{(\ell_2)} |x_\mathsf{v} - y_\mathsf{v}|^2 + s^{(\delta)} \sum_{(\mathsf{v},\mathsf{u})\in\mathsf{E}} \lambda_{\mathsf{uv}}^{(\delta_0)} |x_\mathsf{u} - x_\mathsf{v}|_0 + s^{(\ell)} \sum_{\mathsf{v}\in\mathsf{V}} \lambda_\mathsf{v}^{(\ell_0)} |x_\mathsf{v}|_0 \,,$$

and by analogy with the $\ell_0$-pseudonorm, we denote for any $x \in \mathbb{R}$, $|x|_0 \overset{\text{def}}{=} 1$ if $x \neq 0$ and $|0|_0 = 0$.

| Dataset | | population | revenue | election |
|---|---|---|---|---|
| Observed statistic | $y$ | population density | average revenue | election results |
| Extensive field | $v$ | region surface | population | number of voters |
| Space division | | rasters | rasters | constituencies |
| Spatial extent | | Île-de-France | Île-de-France | France |
| Number of vertices | $\lvert V \rvert$ | 252 183 | 252 183 | 4 670 492 |
| Number of edges | $\lvert E \rvert$ | 378 258 | 378 258 | 7 002 424 |
| Presence of zero extensive quantities | $s^{(\ell)} > 0$ | no | yes | no |

Table 1: Dataset summary for each experimental setting.

The objective functional $F_0$ is the sum of three terms, weighted by their respective penalization coefficients. The first term is a data-fidelity measure, favoring a solution close to the observation. We naturally weight it by the extensive quantity, *i.e.* for all $v \in V$, $\lambda_v^{(\ell_2)}$ is set to $v_v$ if $y_v$ is observed, and to 0 if $y_v$ is missing.

The second term is a penalization ensuring the simplicity of the solution, which tends to merge together neighboring subregions with similar values. We weight the contribution of each pair proportionally to the length of the edge border, *i.e.* for all $(u, v) \in E$, $\lambda_{uv}^{(\delta_0)} \overset{\text{set}}{=} \mu_{uv}$. This term is thus proportional to the total length of the contours of the mapping $x$, in a similar fashion to the geometric term in the *Mumford-Shah functional* (see for instance the review of Vitti, 2012). The coefficient $s^{(\delta)} \in \mathbb{R}_{++}$ scales its influence relatively to the other terms in $F_0$.

Finally, the last term penalizes nonzero values attributed to regions whose extensive quantity is zero. Without this term, large areas could take values of little significance, eventually cluttering the map. Hence, we penalize such regions proportionally to their surface, *i.e.* for all $v \in V$, $\lambda_v^{(\ell_0)}$ is set to $\lambda_v$ if $v_v = 0$, and to 0 otherwise. Again, $s^{(\ell)} \in \mathbb{R}_{++}$ scales its overall influence in $F_0$.

The minimization of $F_0$ is very challenging because the functional $|\cdot|_0$ is noncontinuous and nonconvex. Thus, we consider a *convex relaxation*, defined for all $x \in \mathbb{R}^{|V|}$, as

$$F(x) \overset{\text{def}}{=} \frac{1}{2} \sum_{v \in V} \lambda_v^{(\ell_2)} |x_v - y_v|^2 + \sum_{(v,u) \in E} \lambda_{uv}^{(\delta_1)} |x_u - x_v| + \sum_{v \in V} \lambda_v^{(\ell_1)} |x_v| , \qquad (4.1)$$

where the penalization coefficients $\lambda^{(\delta_1)}$ and $\lambda^{(\ell_1)}$ can be related to $\lambda^{(\delta_0)}$ and $s^{(\delta)}$, and $\lambda^{(\ell_0)}$ and $s^{(\ell)}$, respectively.

## 4.2 Applying a Preconditioned Generalized Forward-Backward

Let $\mathcal{H} \overset{\text{set}}{=} \mathbb{R}^{|V|}$, and $\|\cdot\|$ be the Euclidean norm. There is many ways of casting (4.1) as an instance of (1.3); we describe one of them. An implementation in Matlab, which uses a representation of the graphical data adapted to matricial computations, is available at `https://sites.google.com/site/landrieuloic/home-1`.

### 4.2.1 Splitting

We set the smooth part of the splitting as the first term, *i.e.* for all $x \in \mathbb{R}^{|V|}$,

$$f(x) \overset{\text{set}}{=} \frac{1}{2} \sum_{v \in V} \lambda_v^{(\ell_2)} |x_v - y_v|^2 . \qquad (4.2)$$

This is a quadratic functional whose gradient is immediately $\nabla f \colon x \mapsto (\lambda_{\mathsf{v}}^{(\ell_2)}(x_{\mathsf{v}} - y_{\mathsf{v}}))_{\mathsf{v}\in\mathsf{V}}$; it is easy to see with Remark 2.2 that assumption (h1) is satisfied with $L$ being the diagonal operator whose $\mathsf{v}$-th term is $\lambda_{\mathsf{v}}^{(\ell_2)}$ if $\lambda_{\mathsf{v}}^{(\ell_2)} > 0$, and can be any strictly positive real otherwise.

Now, we must split the remaining terms into simple functionals whose proximity operator is easy to compute. Note that our formulation allows the splitting to avoid the introduction of countless, useless variables in the implementation (see § 1.2 and 1.3), which would considerably increase the computational cost and memory requirements. With $\mathsf{E}_+ \overset{\text{def}}{=} \{(\mathsf{u},\mathsf{v}) \in \mathsf{E} | \lambda_{\mathsf{uv}}^{(\delta_1)} > 0\}$ and $\mathsf{V}_+ \overset{\text{def}}{=} \{\mathsf{v} \in \mathsf{V} | \lambda_{\mathsf{v}}^{(\ell_1)} > 0\}$, we can simply set $n \overset{\text{set}}{=} |\mathsf{E}_+| + |\mathsf{V}_+|$, and consider each term of each sum as a separate functional in the splittng. More precisely, we define for all $(\mathsf{u},\mathsf{v}) \in \mathsf{E}_+$, $\mathcal{H}_{\mathsf{uv}} \overset{\text{set}}{=} \{x \in \mathcal{H} \mid \forall\, \mathsf{w} \in \mathsf{V} \setminus \{\mathsf{u},\mathsf{v}\},\, x_{\mathsf{w}} = 0\}$ and for all $\mathsf{v} \in \mathsf{V}_+$, $\mathcal{H}_{\mathsf{v}} \overset{\text{set}}{=} \{x \in \mathcal{H} | \forall\, \mathsf{w} \in \mathsf{V} \setminus \{\mathsf{v}\},\, x_{\mathsf{w}} = 0\}$. Then, we index the simple functionals by the sets $\mathsf{E}_+$ and $\mathsf{V}_+$ instead of $\{1,\dots,n\}$, where for all $(\mathsf{u},\mathsf{v}) \in \mathsf{E}_+$ and all $\mathsf{w} \in \mathsf{V}_+$,

$$g_{\mathsf{uv}} \colon x \mapsto \lambda_{\mathsf{uv}}^{(\delta_1)} |x_{\mathsf{u}} - x_{\mathsf{v}}| \,, \qquad \text{and} \qquad g_{\mathsf{w}} \colon x \mapsto \lambda_{\mathsf{w}}^{(\ell_1)} |x_{\mathsf{w}}| \,.$$

Note that each $\mathcal{H}_{\mathsf{uv}}$ (respectively each $\mathcal{H}_{\mathsf{v}}$) is isomorphic to $\mathbb{R}^2$ (respectively to $\mathbb{R}$), and that the space of auxiliary variables in Algorithm 1 reduces to $(\times_{(\mathsf{u},\mathsf{v})\in\mathsf{E}_+} \mathcal{H}_{\mathsf{uv}}) \times (\times_{\mathsf{v}\in\mathsf{V}_+} \mathcal{H}_{\mathsf{v}})$, isomorphic to $\mathbb{R}^{2|\mathsf{E}_+|+|\mathsf{V}_+|}$.

The following proposition shows that the proximity operator of each $g_{\mathsf{uv}}$ and of each $g_{\mathsf{v}}$ is easy to compute in any diagonal metric.

**Proposition 4.1.** *Let* $(\lambda,\mu) \in \mathbb{R}_{++}^2$, $g \colon \mathbb{R} \to \mathbb{R} \colon x \mapsto \lambda|x|$, $h \colon \mathbb{R}^2 \to \mathbb{R} \colon (x_1,x_2) \mapsto \mu|x_1 - x_2|$, $m \in \mathbb{R}_{++}$, $(m_1,m_2) \in \mathbb{R}_{++}^2$ *and* $M \overset{\text{def}}{=} \operatorname{diag}(m_1,m_2)$*. Then,*

(i) $$\forall\, x \in \mathbb{R},\ \operatorname{prox}_g^m(x) = \begin{cases} \left(1 - \frac{\lambda/m}{|x|}\right)x & \text{if } |x| > \lambda/m \,, \\ 0 & \text{if } |x| \le \lambda/m \,. \end{cases}$$

(ii) $$\forall \begin{pmatrix} x_1 \\ x_2 \end{pmatrix} \in \mathbb{R}^2,\ \operatorname{prox}_h^M \begin{pmatrix} x_1 \\ x_2 \end{pmatrix} = \begin{cases} \begin{pmatrix} \bar{x} \\ \bar{x} \end{pmatrix} + \left(1 - \frac{\bar{\mu}}{|x_1 - x_2|}\right)\begin{pmatrix} w_2(x_1 - x_2) \\ w_1(x_2 - x_1) \end{pmatrix} & \text{if } |x_1 - x_2| > \bar{\mu} \,, \\ \begin{pmatrix} \bar{x} \\ \bar{x} \end{pmatrix} & \text{if } |x_1 - x_2| \le \bar{\mu} \,, \end{cases}$$

*where* $w_1 \overset{\text{def}}{=} \frac{m_1}{m_1+m_2}$, $w_2 \overset{\text{def}}{=} \frac{m_2}{m_1+m_2}$, $\bar{x} \overset{\text{def}}{=} w_1 x_1 + w_2 x_2$*, and* $\bar{\mu} \overset{\text{def}}{=} \mu(\frac{1}{m_1} + \frac{1}{m_2})$.

*Proof.* (i). Let $x \in \mathbb{R}$. Then, $\operatorname{prox}_g^m(x) = \operatorname{argmin}_{y\in\mathbb{R}} \{\frac{m}{2}(x-y)^2 + \lambda|y|\}$; we do not change the minimizer by dividing the objective by $m$, and we obtain the soft-thresholding.
(ii). Let $x \in \mathbb{R}^2$. If $x_1 = x_2$, then $\bar{x} = x_1 = x_2$ and the result obviously holds. Suppose now $x_1 \ne x_2$. By definition, $\operatorname{prox}_h^M(x) = \operatorname{argmin}_{y\in\mathbb{R}^2}\{H(y)\}$ where for all $y \in \mathbb{R}^2$, $H(y) \overset{\text{def}}{=} \frac{m_1}{2}(x_1 - y_1)^2 + \frac{m_2}{2}(x_2 - y_2)^2 + \mu|y_1 - y_2|$. For all $y \in \mathbb{R}^2$ such that $y_1 \ne y_2$, $H$ is differentiable in $y$, with

$$\frac{\partial H}{\partial y_1}(y) = m_1(y_1 - x_1) + \mu\frac{y_1 - y_2}{|y_1 - y_2|} \,, \quad \text{and} \quad \frac{\partial H}{\partial y_2}(y) = m_2(y_2 - x_2) + \mu\frac{y_2 - y_1}{|y_1 - y_2|} \,.$$

We deduce from first-order optimality conditions that $y \in \mathbb{R}^2$ such that $y_1 \ne y_2$ is equal to $\operatorname{prox}_h^M(x)$ if, and only if,

$$y_1 = x_1 - \frac{\mu}{m_1}\frac{y_1 - y_2}{|y_1 - y_2|} \,, \quad \text{and} \quad y_2 = x_2 - \frac{\mu}{m_2}\frac{y_2 - y_1}{|y_1 - y_2|} \,. \tag{4.3}$$

Supposing that such a $y$ exists, it necessarily satisfies $y_1 - y_2 = x_1 - x_2 - \bar{\mu}\frac{y_1-y_2}{|y_1-y_2|}$. That is, if $y_1 > y_2$, then $x_1 - x_2 > \bar{\mu}$, and if $y_1 < y_2$, then $x_1 - x_2 < -\bar{\mu}$. Altogether, a necessary condition is $|x_1 - x_2| > \bar{\mu}$. But in that case, observe that with $y_1 \overset{\text{def}}{=} x_1 - \frac{\mu}{m_1}\frac{x_1-x_2}{|x_1-x_2|}$ and $y_2 \overset{\text{def}}{=} x_2 - \frac{\mu}{m_2}\frac{x_2-x_1}{|x_1-x_2|}$, $y_1 - y_2$ has the same sign as $x_1 - x_2$, and thus $y \overset{\text{def}}{=} (y_1, y_2)$ satisfies (4.3). In order to obtain the final expression, note for instance that $x_1 - \bar{x} = w_2(x_1 - x_2)$, so that $x_1 - \frac{\mu}{m_1}\frac{x_1-x_2}{|x_1-x_2|} = \bar{x} + w_2(x_1 - x_2)\big(1 - \frac{\mu(m_1+m_2)}{|x_1-x_2|m_1 m_2}\big)$ and the result follows.
If at the contrary $|x_1 - x_2| \leq \bar{\mu}$, we deduce by contraposition that $\text{prox}_h^M(x) = (\bar{y}, \bar{y})$, where $\bar{y} \overset{\text{def}}{=} \text{argmin}_{y\in\mathbb{R}}\{\frac{m_1}{2}(x_1 - y)^2 + \frac{m_2}{2}(x_2 - y)^2\}$, and optimality conditions immediately lead to $\bar{y} = \bar{x}$. ∎

Finally, note that we can assume without loss of generality that $\text{span}\big\{\{\mathcal{H}_{\mathsf{uv}}\,|(\mathsf{u},\mathsf{v}) \in \mathsf{E}_+\} \cup \{\mathcal{H}_{\mathsf{v}}\,|\mathsf{v} \in \mathsf{V}_+\}\big\} = \mathcal{H}$. Otherwise, this would mean that there exists a vertex $\mathsf{v} \in \mathsf{V} \setminus \mathsf{V}_+$ such that for all $\mathsf{u} \in \mathsf{V}$, $(\mathsf{u},\mathsf{v}) \notin \mathsf{E}_+$ and $(\mathsf{v},\mathsf{u}) \notin \mathsf{E}_+$. But then the label $x_{\mathsf{v}}$ of such a vertex intervenes in the functional $F$ only through the term $\lambda_{\mathsf{v}}^{(\ell_2)}|x_{\mathsf{v}} - y_{\mathsf{v}}|^2$, in which case it can be set directly to $y_{\mathsf{v}}$ if it is observed, or simply discarded from the problem otherwise.

### 4.2.2 Preconditioning

We follow here § 3.2. Let us first obtain suitable quadratic approximations. As already mentioned, $f$ is itself quadratic with a diagonal Hessian. The latter is not necessarily strongly positive, but as explained below, this is not a concern.

Let now $\mathsf{v} \in \mathsf{V}_+$, and recall that $g_{\mathsf{v}}\colon x \mapsto \lambda_{\mathsf{v}}^{(\ell_1)}|x_{\mathsf{v}}|$. Given $\hat{x} \in \mathcal{H}$, we consider the approximation of $g_{\mathsf{v}}$ at $\hat{x}$ as a quadratic function $q_{\mathsf{v}}$ such that $q_{\mathsf{v}}(\hat{x}) = g_{\mathsf{v}}(\hat{x})$, $\nabla q_{\mathsf{v}}(\hat{x}) \in \partial g_{\mathsf{v}}(\hat{x})$, and for all $x \in \mathcal{H}$, $q_{\mathsf{v}}(x) \geq g_{\mathsf{v}}(x)$. As long as $\hat{x}_{\mathsf{v}} \neq 0$, it is easy to see that the best candidate is given by

$$q_{\mathsf{v}}\colon x \mapsto \frac{\lambda_{\mathsf{v}}^{(\ell_1)}}{2}\frac{{x_{\mathsf{v}}}^2}{|\hat{x}_{\mathsf{v}}|} + \frac{\lambda_{\mathsf{v}}^{(\ell_1)}}{2}|\hat{x}_{\mathsf{v}}|\,,$$

as illustrated on Figure 1. Its Hessian is the diagonal operator whose only nonzero term is the $\mathsf{v}$-th, equal to $\lambda_{\mathsf{v}}^{(\ell_1)}/|\hat{x}_{\mathsf{v}}|$; its restriction to $\mathcal{H}_{\mathsf{v}}$ belongs to $\mathfrak{S}_{++}(\mathcal{H}_{\mathsf{v}})$. However, the latter term explodes as $\hat{x}_{\mathsf{v}}$ tends to 0, inducing instability. To ensure robustness of the reconditioning procedure, we replace $|\hat{x}_{\mathsf{v}}|$ with $\max\{|\hat{x}_{\mathsf{v}}|, \epsilon_{\ell_1}\}$, for some small value $\epsilon_{\ell_1} \in \mathbb{R}_{++}$. In all experiments, we observed good robustness by setting it as a billionth of the average amplitude of the values, *i.e.* $\epsilon_{\ell_1} \overset{\text{set}}{=} \frac{10^{-6}}{|\mathsf{V}|}\sum_{\mathsf{v}\in\mathsf{V}}|\hat{x}_{\mathsf{v}}|$.

Letting now $(\mathsf{u},\mathsf{v}) \in \mathsf{E}_+$, by analogy with the absolute value above, as long as $\hat{x}_{\mathsf{u}} \neq \hat{x}_{\mathsf{v}}$ we consider the quadratic approximation of $g_{\mathsf{uv}}$ at $\hat{x}$ as

$$q_{\mathsf{uv}}\colon x \mapsto \frac{\lambda_{\mathsf{uv}}^{(\delta_1)}}{2}\frac{(x_{\mathsf{u}} - x_{\mathsf{v}})^2}{|\hat{x}_{\mathsf{u}} - \hat{x}_{\mathsf{v}}|} + \frac{\lambda_{\mathsf{uv}}^{(\delta_1)}}{2}|\hat{x}_{\mathsf{u}} - \hat{x}_{\mathsf{v}}|\,.$$

In contrast however, its Hessian is not diagonal. Recalling § 3.2.1, we apply the infinite-norm equilibration by simply dropping the off-diagonal terms. The result is the diagonal operator whose only nonzero terms are the $\mathsf{u}$-th and $\mathsf{v}$-th, both equal to $\lambda_{\mathsf{uv}}^{(\delta_1)}/|\hat{x}_{\mathsf{u}} - \hat{x}_{\mathsf{v}}|$; its restriction to $\mathcal{H}_{\mathsf{uv}}$ belongs to $\mathfrak{S}_{++}(\mathcal{H}_{\mathsf{uv}})$. Once again the latter term explodes as $\hat{x}_{\mathsf{u}} - \hat{x}_{\mathsf{v}}$ tends to 0, and we must replace $|\hat{x}_{\mathsf{u}} - \hat{x}_{\mathsf{v}}|$ with $\max\{|\hat{x}_{\mathsf{u}} - \hat{x}_{\mathsf{v}}|, \epsilon_{\delta_1}\}$ for some small value $\epsilon_{\delta_1} \in \mathbb{R}_{++}$. In contrast with $\epsilon_{\ell_1}$, we observed low robustness with respect to $\epsilon_{\delta_1}$ on several experimental settings. In order to avoid parameter tweaking, we set it for all experiments to a fairly large value, namely $\epsilon_{\delta_1} \stackrel{\text{set}}{=} \max\{|\hat{x}_u|/10, \epsilon_{\ell_1}\}$.

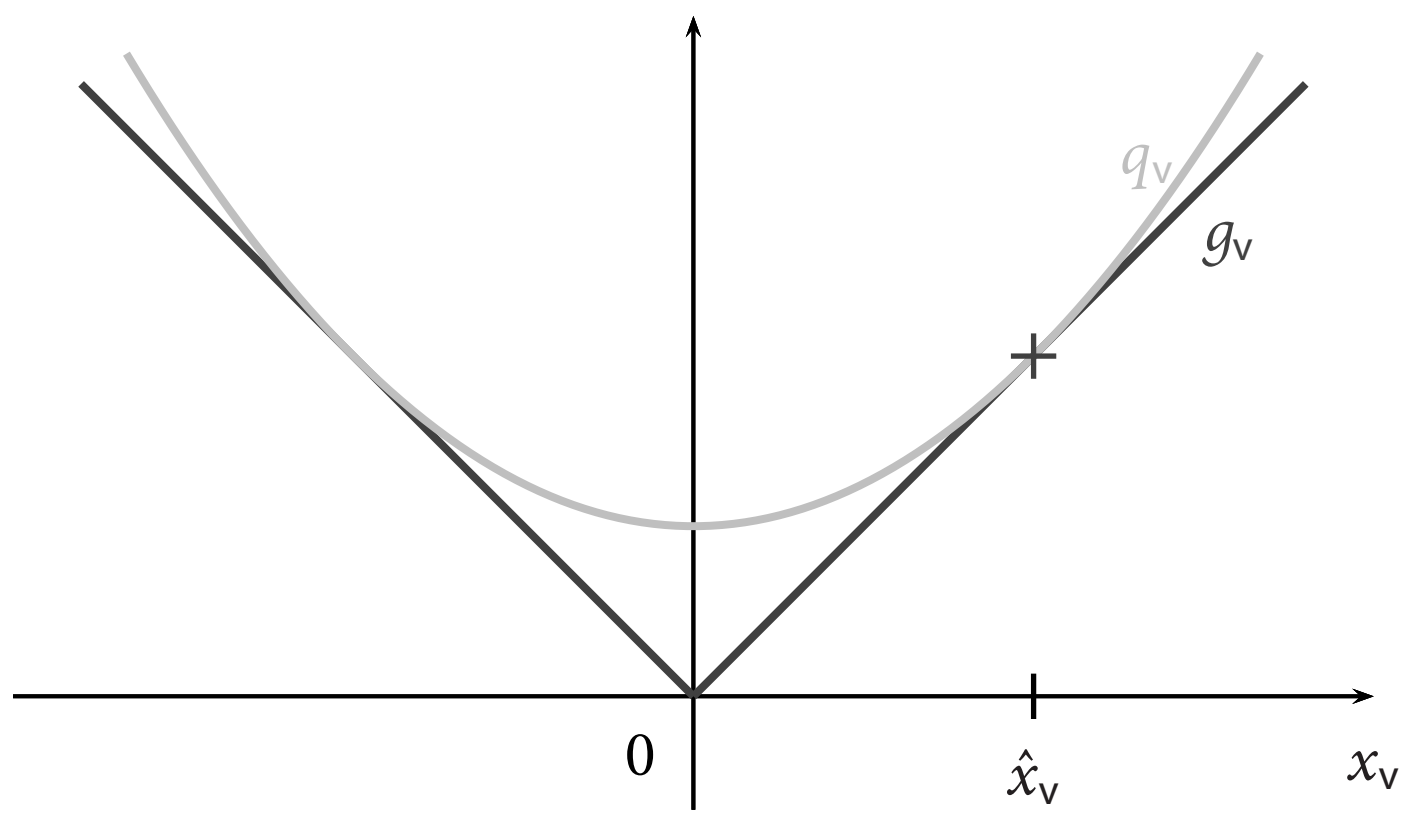

Figure 1: Quadratic approximation of $g_v$.

Altogether, we are now set to apply the preconditioning procedure explained along § 3.2. In particular, we set for all $k \in \mathbb{N}$, $\rho_k \stackrel{\text{set}}{=} 1.5$, and $\Gamma$ according to (3.1b); since we have span $\{\{\mathcal{H}_{uv} \mid (u,v) \in E_+\} \cup \{\mathcal{H}_v \mid v \in V_+\}\} = \mathcal{H}$, this ensures in addition $\Gamma \in \mathfrak{S}_{++}(\mathcal{H})$. Finally, since with Proposition 4.1 the needed proximity operators can be computed in any diagonal metric, we set $(W_v)_{v \in V_+}$ and $(W_{uv})_{(u,v) \in E_+}$ according to (3.3a).

## 4.3 Experimental Setup and Results

We perform aggregation of spatial statistics over three different datasets, presented in Table 1. The datasets `population` and `revenue` are open-source, available at the French national institute for statistic and economics studies[1]. The dataset `election` is also open-source, provided by the *Cartelec* project (Colange *et al.*, 2013).

The spatial regions of rasterized data are squares organized along regular grids, while the constituencies in `election` constitute arbitrary polygons. More readable maps are obtained by Delaunay triangulation on all vertices composing the shape of regions, with the constraint that all region borders must be used as edges by the triangles (see Chew, 1989). Each region is thus split into triangles, to which we attribute the intensive value observed on the region. In turn, the extensive values is shared among the triangles, proportionally to their surface.

In the following, we present the performance of our preconditioned generalized forward-backward splitting algorithm over minimizations of $F$ as defined in (4.1), against other preconditioned proximal splitting algorithms available in the literature. For the sake of completeness and illustration, we then briefly comment some results obtained at the end of the overall procedure described all along § 4.1.

### 4.3.1 Optimization Comparison

The following preconditioned proximal splitting algorithms are considered in the comparison. Each algorithm is implemented in MATLAB, and made available at `https://www.ceremade.`

[1] *IdeesLibres.org* 01/2015, INSEE 20/11/2013,
`https://www.data.gouv.fr/fr/datasets/donnees-carroyees-a-200m-sur-la-population/`

dauphine.fr/~raguet/pgfb/, together with a convenient matricial representation of the dataset revenue for reproducibility.

**Preconditioned primal-dual algorithm of** Pock and Chambolle (2011, PPD). The functional is split as $F \overset{\text{def}}{=} f + g \circ K$, where $f$ is the data-fidelity term as in (4.2), and $g$ and $K$ are defined as

$$K\colon \begin{array}{ccc} \mathbb{R}^{|\mathsf{V}|} & \longrightarrow & \mathbb{R}^{|\mathsf{E}|} \times \mathbb{R}^{|\mathsf{V}|} \\ x & \longmapsto & (\delta, \xi)\,, \end{array} \quad \text{with} \quad \left\{ \begin{array}{ll} \forall\,(\mathsf{u},\mathsf{v}) \in \mathsf{E}, & \delta_{\mathsf{uv}} \overset{\text{def}}{=} \lambda^{(\delta_1)}_{\mathsf{uv}} (x_{\mathsf{u}} - x_{\mathsf{v}})\,, \\ \forall\,\mathsf{v} \in \mathsf{V}, & \xi_{\mathsf{v}} \overset{\text{def}}{=} \lambda^{(\ell_1)}_{\mathsf{v}} x_{\mathsf{v}}\,, \end{array} \right.$$

$$\text{and} \quad g\colon \begin{array}{ccc} \mathbb{R}^{|\mathsf{E}|} \times \mathbb{R}^{|\mathsf{V}|} & \longrightarrow & \mathbb{R} \\ (\delta, \xi) & \longmapsto & \sum_{(\mathsf{u},\mathsf{v})\in \mathsf{V}} |\delta_{\mathsf{uv}}| + \sum_{\mathsf{v}\in\mathsf{V}} |\xi_{\mathsf{v}}|\,. \end{array}$$

The proximity operator of the functionals $f$ and $g$ are easy to compute, and $K$ is a linear operator, and thus iteration (4) of Pock and Chambolle (2011) can be applied. The preconditioning matrices $T$ and $\Sigma$ are defined following Lemma 2, equation (10), with the parameter $\alpha \overset{\text{set}}{=} 1$. We also set the relaxation parameter $\theta \overset{\text{set}}{=} 1$.

Note that the preconditioning procedure only depends on the operator $K$; in particular, no information on the functionals $f$ or $g$ is taken into account. Thus, another classical splitting, which consists in taking into account the coefficients $(\lambda^{(\delta_1)}_{\mathsf{uv}})_{(\mathsf{u},\mathsf{v})\in\mathsf{E}}$ and $(\lambda^{(\ell_1)}_{\mathsf{v}})_{\mathsf{v}\in\mathsf{V}}$ within the functional $g$, gives bad performances on our data, where the coefficients present highly heterogeneous values (data not shown). Conversely, we also tried another splitting where the role of $f$ is captured by $K$ and $g$, but no improvement was observed (data not shown either).

**Inertial preconditioned primal-dual algorithm of** Lorenz and Pock (2015, IPPD). The iteration (30) of Lorenz and Pock (2015) can be seen as an inertial extension of the above, where in addition the functional $f$ can be taken into account through an explicit gradient step. The preconditioning matrices $T$ and $\Sigma$ can in turn incorporate information over $f$ following Lemma 10, equation (35), recalling in their terms that $\nabla f$ is cocoercive with respect to the diagonal matrix $L^{-1}$. We used the parameters $\gamma \overset{\text{set}}{=} 1$, $\delta \overset{\text{set}}{=} 0$, $r \overset{\text{set}}{=} 1$, $s \overset{\text{set}}{=} 1$, and, after trying different ratios, the inertial parameters $(\alpha_k)_{k\in\mathbb{N}}$ all equal to one half of the upper bound given by Lemma 6, equation (25).

**Preconditioned generalized forward-backward splitting** (Algorithm 1, PGFB$_\theta$). We use the splitting and preconditioning described along § 4.2. Note that at the very beginning of the optimization process, we initialize the conditioners with a coarse preconditioning, following § 4.2.2 but where each amplitude $|x_{\mathsf{v}}|$ and each absolute difference $|x_{\mathsf{u}} - x_{\mathsf{v}}|$ are replaced by the average amplitude of the observed values in $y$.

Then, we consider reconditionings taking place when the relative evolution of the iterate at iteration $k$, *i.e.* $\|x_k - x_{k-1}\|/\|x_{k-1}\|$, is below a certain threshold. Within the run of the algorithm, when this threshold is reached, it is then divided by 10. We denote the resulting method with given *initial* threshold $\theta$ by PGFB$_\theta$; in particular, PGFB$_0$ means that no reconditioning is performed.

For each dataset, we fix reasonable values for the parameters $s^{(\delta)}$ and $s^{(\ell)}$ (corresponding ultimately to some of the illustrations given in § 4.3.2). Then, for each algorithm, we monitor first the computational time needed for completing a thousand iterations, and then the decrease of the objective functional $F$ along these iterations. Finally, an approximate minimum of $F$ is obtained with five thousands iterations of PGFB. We plot on Figure 2 the distance between the values of the objective along the iterations and this minimum on a logarithmic scale, against the corresponding computational time.

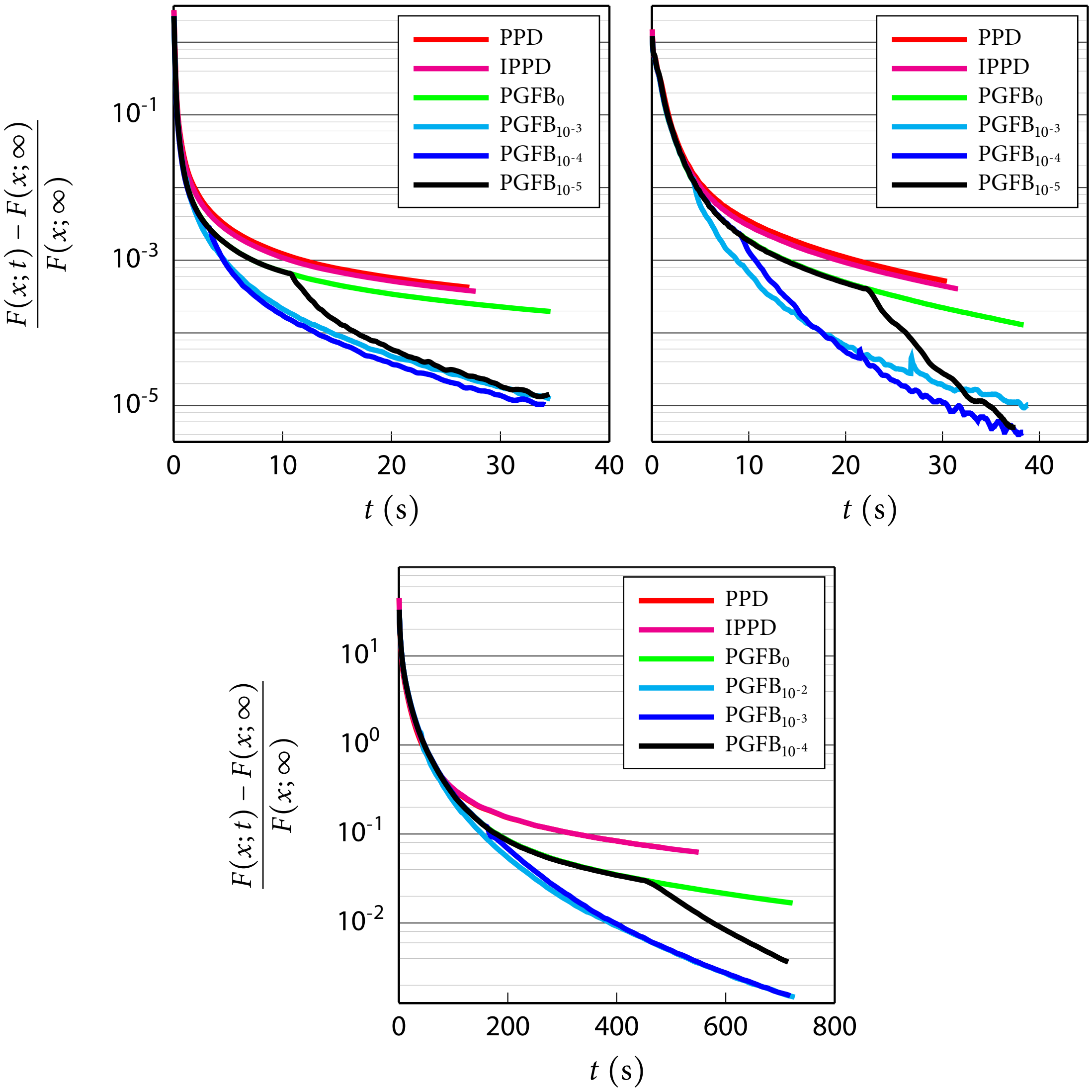

Figure 2: Optimization comparisons. The top row shows the performance over the smallest datasets `population` and `revenue`, while the last graphic shows the performance over the significantly larger dataset `election`.

All performance results show the same tendency, in spite of the variety of the problems, of the conditionings, and of the data size considered. In all three experiments, we see that PPD and IPPD iterations are faster than PGFB iterations; yet, the coarse initial preconditioning is already enough for $\mathrm{PGFB}_0$ to outperform PPD and IPPD.

The three different reconditioning threshold values illustrate well the compromise which must be found for optimal reconditioning—not too early, not too late—evoked at the end of § 3.2.2. However in all cases, it is clear that the computational cost of the reconditionings are negligible, and that they allow to drastically outperform $\mathrm{PGFB}_0$.

Interestingly, we can see on Figure 2, top right (`revenue` dataset), an undesirable jump in the functional values, occurring for the cyan ($\text{PGFB}_{10^{-3}}$) and blue ($\text{PGFB}_{10^{-4}}$) curves between 20 and 30 s. These are due to reconditionings, for reasons explained in § 3.2.2; fortunately, more advantageous functional values are recovered within only a few iterations. Note that the subsequent small fluctuations in the functional values are a known side-effect of the use of overrelaxation ($\rho_k > 1$), and should not be confused with reconditionings.

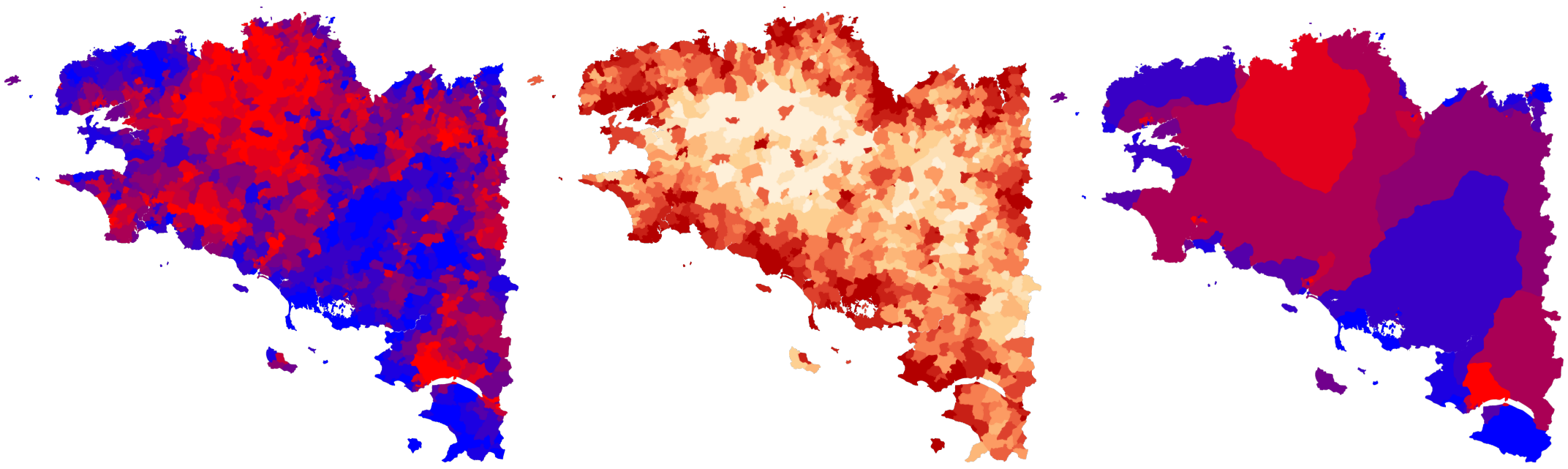

(a) (b) (c)

Figure 3: Close-up on the `election` dataset: results of the second round of the 2007 French presidential election, broken down by constituencies in French Brittany. Two candidates are opposed, the colormap goes thus from blue to red, representing respectively the regions where one candidate achieve its highest score, respectively its lowest score. (a). Original map. (b). Number of voters per surface unit over constituencies, from low density in pale orange to high density in dark red. (c). Aggregation; compression ratio 9, relative error 0.22.

### 4.3.2 Interpretation of the Results

We illustrate and comment briefly on the usefulness of our approach for aggregation of spatial statistics. More substantive discussion should be found elsewhere (Landrieu *et al.*, 2015). Let us precise here that solving only one convex relaxation of $F_0$ often leads to undesirable *staircasing* and *bias* effects. To reduce them, we actually solve successively several instances of $F$ as defined in (4.1), but with coefficients $\lambda^{(\delta_1)}$ and $\lambda^{(\ell_1)}$ depending each time on the previously found solution, following classical *reweighting* techniques (see in particular the recent review of Ochs *et al.*, 2015).

As long as there is no sparsity-inducing $\ell_1$ penalization, our aggregation method can be seen as a lossy compression process, where one simply seeks for a trade-off between data size and loss of information. As already pointed out in the aggregation model, an interesting measure of the complexity of a map is given by the total length of the contours between regions. We thus measure the *compression ratio* of the aggregation $x$ of the spatial statistics $y$ as

$$\frac{\sum_{(\mathsf{v},\mathsf{u})\in E} \mu_{\mathsf{uv}} |y_{\mathsf{u}} - y_{\mathsf{v}}|_0}{\sum_{(\mathsf{v},\mathsf{u})\in E} \mu_{\mathsf{uv}} |x_{\mathsf{u}} - x_{\mathsf{v}}|_0} .$$

At the same time, a relevant measure of the *relative error* is the root (weighted) mean square error between simplified and observed map, which we normalize by the (weighted) standard-deviation of the latter, *i.e.*

$$\frac{\sqrt{\sum_{\mathsf{v}\in V} \nu_{\mathsf{v}} (x_{\mathsf{v}} - y_{\mathsf{v}})^2}}{\sqrt{\sum_{\mathsf{v}\in V} \nu_{\mathsf{v}} (y_{\mathsf{v}} - \bar{y})^2}} , \quad \text{where} \quad \bar{y} \overset{\text{def}}{=} \frac{\sum_{\mathsf{v}\in V} \nu_{\mathsf{v}} y_{\mathsf{v}}}{\sum_{\mathsf{v}\in V} \nu_{\mathsf{v}}} .$$

Such measures are reported on Figures 3 and 4. Note in particular on the former that our model yields final region sizes that adapt to the local density, overlook differences between low density area in favors of statistically more significant local effect such as political polarization at the scale of a city.

Because of the presence of the sparsity-inducing penalization term, the aggregations on the `revenue` dataset are somewhat more difficult to interpret. Note that on Figure 5 (b) to (d), local areas without population can still be distinguished, in spite of high degrees of simplification. Then, as can be seen on Figure 5 (e) and (f), even larger penalizations can be used to suppress scarcely distributed information instead of averaging it.

# 5 Convergence Proof

The convergence proof follows the same line as our generalized forward-backward (Raguet *et al.*, 2013, Theorem 2.1), with the appropriate modifications. We begin by recalling important properties in monotone operators theory. We give useful definitions and properties over the product space $\times_{i=1}^{n} \mathcal{H}_i$. Finally, we formulate a fixed-point equation satisfied by the solutions of (1.1), from which we derive the algorithmic scheme and its convergence.

## 5.1 Preliminaries on Monotone Operators Theory

The following lemma is useful for manipulating averaged and cocoercive operators.

**Lemma 5.1.** *Let* $(\alpha, \alpha') \in \left]0,1\right[^2$, $\beta \in \left]0,+\infty\right[$, $\gamma \in \left]0,2\beta\right[$, *and* $T, T' \colon \mathcal{H} \to \mathcal{H}$.

(i) *$T$ is $\alpha$-averaged in $\mathcal{H}$, if, and only if, for all $(x,y) \in \mathcal{H}^2$,*

$$\|Tx - Ty\|^2 \leq \|x-y\|^2 - \tfrac{1-\alpha}{\alpha}\|(\mathrm{Id}-T)x - (\mathrm{Id}-T)y\|^2\,;$$

(ii) *$T$ is $\alpha$-averaged in $\mathcal{H}$, if, and only if, for all $(x,y) \in \mathcal{H}^2$,*

$$2(1-\alpha)\langle Tx - Ty \,|\, x-y\rangle \geq \|Tx-Ty\|^2 + (1-2\alpha)\|x-y\|^2\,;$$

(iii) *If $T \in \mathcal{A}(\alpha,\mathcal{H})$ and $T' \in \mathcal{A}(\alpha',\mathcal{H})$, then $TT' \in \mathcal{A}(\frac{\alpha+\alpha'-2\alpha\alpha'}{1-\alpha\alpha'},\mathcal{H})$;*

(iv) *$T$ is $\beta$-cocoercive in $\mathcal{H}$ if, and only if, $\mathrm{Id}-\gamma T$ is $\frac{\gamma}{2\beta}$-averaged in $\mathcal{H}$.*

*Proof.* (i)–(ii). See Combettes (2004, Lemma 2.1 (i) ⇔ (ii) (i) ⇔ (iii)).
(iii). See Ogura and Yamada (2002, Theorem 3). These authors actually allow $\alpha \in \left[0,1\right[$; the above result is derived from theirs with $\alpha + \alpha' \geq 2\min\{\alpha,\alpha'\} > 2\alpha\alpha'$, ensuring $\frac{\alpha+\alpha'-2\alpha\alpha'}{1-\alpha\alpha'} > 0$.
(iv). $\beta T = \frac{1}{2}(R+\mathrm{Id}) \Leftrightarrow \mathrm{Id}-\gamma T = \mathrm{Id}-\frac{\gamma}{2\beta}(R+\mathrm{Id}) = \frac{\gamma}{2\beta}(-R) + \left(1-\frac{\gamma}{2\beta}\right)\mathrm{Id}$. ∎

The following lemma is useful for characterizing firmly nonexpansive operators.

**Lemma 5.2.** *The following statements are equivalent:*

(i) *$T$ is firmly nonexpansive in $\mathcal{H}$ with full domain (i.e.* $\operatorname{dom} T = \mathcal{H}$*);*

(ii) *$2T - \mathrm{Id}$ is nonexpansive in $\mathcal{H}$ and $T$ has full domain;*

(iii) *$\forall\, (x,y) \in \mathcal{H}^2$, $\|Tx-Ty\|^2 \leq \langle Tx - Ty \,|\, x-y\rangle$;*

(iv) *$T$ is the resolvent of a maximally monotone operator $A$, i.e. $T = J_A$.*

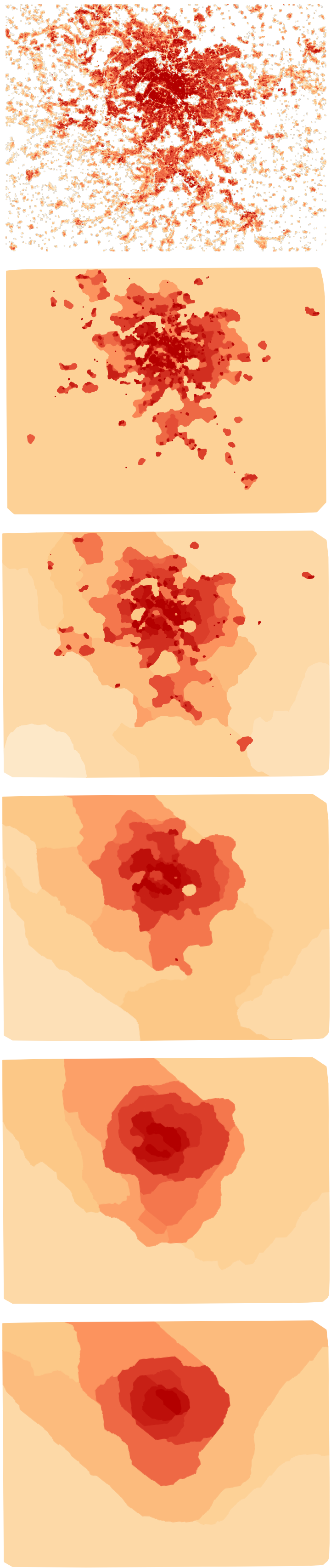

Figure 4

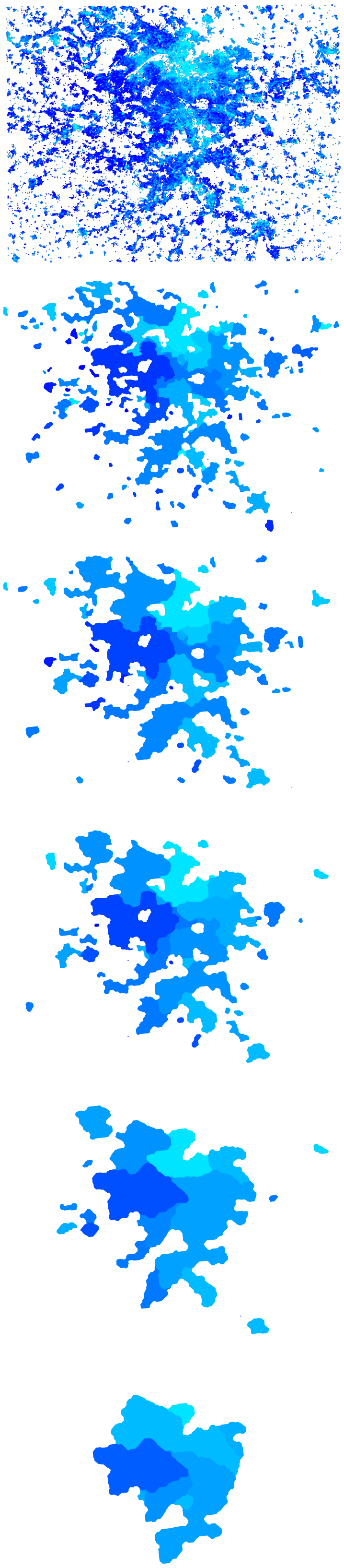

Figure 5

(a)
(b)
(c)
(d)
(e)
(f)

(a). Original maps. (b)–(f). Aggregations with increasing degrees of simplification.

Figure 4: Aggregation of the population density in the greater Paris area for increasing values of $s^{(\delta)}$. The colormap represents high density areas in dark red and low density areas in pale orange. The table below presents the evolution of the compression and error.

| | comp. | err. |
|---|---|---|
| (a) | 1 | 0 |
| (b) | 6 | 0.34 |
| (c) | 12 | 0.45 |
| (d) | 21 | 0.51 |
| (e) | 26 | 0.54 |
| (f) | 38 | 0.57 |

Figure 5: Aggregation of the average yearly revenue density in the greater Paris area for increasing values of $s^{(\delta)}$ and $s^{(\ell)}$. The colormap represents areas of high revenues in dark blue and areas of low revenues in cyan.

*Proof.* (i)⇔(ii). $T \in \mathcal{A}\left(\frac{1}{2}, \mathcal{H}\right) \Leftrightarrow T = \frac{1}{2}(R + \mathrm{Id})$ for some $R$ nonexpansive. (ii)⇔(iii). This is Lemma 5.1 (ii) with $\alpha \overset{\text{set}}{=} \frac{1}{2}$, introduced in this particular case by Zarantonello (1971, Lemma 1.3). (i)⇔(iv). See Minty (1962) and again Zarantonello (1971). ■

## 5.2 Definitions and Properties on the Product Space

Let $\boldsymbol{\mathcal{H}} \overset{\text{def}}{=} \times_{i=1}^{n} \mathcal{H}_i$, whose generic element is denoted $\boldsymbol{x} \overset{\text{def}}{=} (x_i)_{1 \leq i \leq n}$, endowed with the canonical inner product defined for all $(\boldsymbol{x}, \boldsymbol{y}) \in \boldsymbol{\mathcal{H}}^2$ by $\langle \boldsymbol{x} \,|\, \boldsymbol{y} \rangle \overset{\text{def}}{=} \sum_{i=1}^{n} \langle x_i \,|\, y_i \rangle$, and induced norm $\|\boldsymbol{x}\| \overset{\text{def}}{=} \sqrt{\langle \boldsymbol{x} \,|\, \boldsymbol{x} \rangle} = \sqrt{\sum_{i=1}^{n} \|x_i\|^2}$. We particularize $\mathbf{Id} \overset{\text{def}}{=} \mathrm{Id}_{\boldsymbol{\mathcal{H}}}$, and with the definitions from § 2.2, we define the following linear operators $\boldsymbol{\mathcal{H}} \to \boldsymbol{\mathcal{H}}$, $\boldsymbol{W} \colon \boldsymbol{x} \mapsto (W_i x_i)_{1 \leq i \leq n}$ and $\boldsymbol{\Gamma} \colon \boldsymbol{x} \mapsto (\Gamma x_i)_{1 \leq i \leq n}$, whose respective inverses are $\boldsymbol{W}^{-1} \colon \boldsymbol{x} \mapsto (W_i^{-1} x_i)_{1 \leq i \leq n}$ and $\boldsymbol{\Gamma}^{-1} \colon \boldsymbol{x} \mapsto (\Gamma^{-1} x_i)_{1 \leq i \leq n}$. Moreover, we also introduce

$$\boldsymbol{\Sigma} \colon \boldsymbol{\mathcal{H}} \to \mathcal{H} \colon \boldsymbol{x} \mapsto \sum_{i=1}^{n} W_i x_i \qquad \text{and} \qquad \Sigma^* \colon \mathcal{H} \to \boldsymbol{\mathcal{H}} \colon x \mapsto \left(x^{\mathcal{H}_i}\right)_{1 \leq i \leq n},$$

together with the set $\boldsymbol{\mathcal{S}} \overset{\text{def}}{=} \operatorname{ran} \Sigma^* = \left\{ \boldsymbol{x} \in \boldsymbol{\mathcal{H}} \,\middle|\, \exists x \in \mathcal{H}, \forall\, i \in \{1, \ldots, n\}, x_i = x^{\mathcal{H}_i} \right\}$. Finally, we define $\boldsymbol{A} \colon \boldsymbol{\mathcal{H}} \to 2^{\boldsymbol{\mathcal{H}}}$ and $\boldsymbol{B} \colon \boldsymbol{\mathcal{H}} \to \boldsymbol{\mathcal{H}}$ by

$$\boldsymbol{A} \overset{\text{def}}{=} \bigtimes_{i=1}^{n} A_{i\,|\mathcal{H}_i} \qquad \text{and} \qquad \boldsymbol{B} \overset{\text{def}}{=} \Sigma^* B \boldsymbol{\Sigma}\,;$$

that is, the graph of $\boldsymbol{A}$ is $\operatorname{gra} \boldsymbol{A} \overset{\text{def}}{=} \times_{i=1}^{n} \operatorname{gra} A_{i|\mathcal{H}_i} = \left\{ (\boldsymbol{x}, \boldsymbol{y}) \in \boldsymbol{\mathcal{H}}^2 \,\middle|\, \forall\, i \in \{1, \ldots, n\}, y_i \in A_i x_i \right\}$. The following proposition sums up the main properties of the above operators.

**Lemma 5.3.** *With the above definition, let* $\boldsymbol{s} \in \boldsymbol{\mathcal{S}}$, $\boldsymbol{x} \in \boldsymbol{\mathcal{H}}$, *and* $x \in \mathcal{H}$. *Then,*

(i) $\boldsymbol{\Sigma} \Sigma^* x = x$,

(ii) $\Sigma^* \boldsymbol{\Sigma} \boldsymbol{s} = \boldsymbol{s}$

(iii) $\boldsymbol{\Gamma} \boldsymbol{s} = \Sigma^* \Gamma \boldsymbol{\Sigma} \boldsymbol{s}$.

(iv) $\langle \boldsymbol{s} \,|\, \boldsymbol{x} \rangle_{\boldsymbol{\Gamma}^{-1}\boldsymbol{W}} = \langle \Gamma^{-1} \boldsymbol{\Sigma} \boldsymbol{s} \,|\, \boldsymbol{\Sigma} \boldsymbol{x} \rangle$.

*Proof.* (i). $\boldsymbol{\Sigma} \Sigma^* x \overset{\text{def}}{=} \sum_{i=1}^{n} W_i x^{\mathcal{H}_i} = \sum_{i=1}^{n} W_i x = x$, where we used (P2) (i) and (iv).
(ii). By definition, there exists $s \in \mathcal{H}$ such that $\boldsymbol{s} = \Sigma^* s$. Using (i), $\boldsymbol{\Sigma} \boldsymbol{s} = s$, and the result follows.
(iii). With (P1) (ii), we have $\boldsymbol{\Gamma} \Sigma^* = \Sigma^* \Gamma$, and thus using (ii), $\boldsymbol{\Gamma} \boldsymbol{s} = \boldsymbol{\Gamma} \Sigma^* \boldsymbol{\Sigma} \boldsymbol{s} = \Sigma^* \Gamma \boldsymbol{\Sigma} \boldsymbol{s}$.
(iv). Take $s = \boldsymbol{\Sigma} \boldsymbol{s}$ as above, and using the properties (P1)–(P2), develop

$$\begin{aligned} \langle \boldsymbol{s} \,|\, \boldsymbol{x} \rangle_{\boldsymbol{\Gamma}^{-1}\boldsymbol{W}} &= \langle \boldsymbol{\Gamma}^{-1} \boldsymbol{W} \Sigma^* s \,|\, \boldsymbol{x} \rangle = \sum_{i=1}^{n} \langle \Gamma^{-1} W_i s^{\mathcal{H}_i} \,|\, x_i \rangle = \sum_{i=1}^{n} \langle \Gamma^{-1} W_i s \,|\, x_i \rangle\,, \\ &= \sum_{i=1}^{n} \langle W_i \Gamma^{-1} s \,|\, x_i \rangle = \sum_{i=1}^{n} \langle \Gamma^{-1} s \,|\, W_i x_i \rangle = \langle \Gamma^{-1} s \,|\, \textstyle\sum_{i=1}^{n} W_i x_i \rangle = \langle \Gamma^{-1} \boldsymbol{\Sigma} \boldsymbol{s} \,|\, \boldsymbol{\Sigma} \boldsymbol{x} \rangle\,. \end{aligned}$$ ■

**Proposition 5.1.** *With the above definitions, the following statements hold.*

(i) $\boldsymbol{\Gamma} \in \mathfrak{S}_{++}(\boldsymbol{\mathcal{H}})$, $\boldsymbol{W} \in \mathfrak{S}_{++}(\boldsymbol{\mathcal{H}})$, *and* $\boldsymbol{\Gamma}^{-1} \boldsymbol{W} \in \mathfrak{S}_{++}(\boldsymbol{\mathcal{H}})$.

(ii) $\boldsymbol{\mathcal{S}}$ *is a closed subspace of* $\boldsymbol{\mathcal{H}}_{\boldsymbol{\Gamma}^{-1}\boldsymbol{W}}$, *and* $\Sigma^* \boldsymbol{\Sigma} = P_{\boldsymbol{\mathcal{S}}}^{\boldsymbol{\Gamma}^{-1}\boldsymbol{W}}$.

(iii) $J_{\boldsymbol{W}^{-1}\boldsymbol{\Gamma A}} = \times_{i=1}^{n} J_{W_i^{-1}\Gamma A_i|\mathcal{H}_i}$, *and* $\boldsymbol{W}^{-1}\boldsymbol{\Gamma A}$ *is maximally monotone in* $\boldsymbol{\mathcal{H}}_{\boldsymbol{\Gamma}^{-1}\boldsymbol{W}}$.

(iv) *If* $B$ *is uniformly monotone in* $\mathcal{H}$ *with modulus* $\phi$, *then*

$$\forall\,(\boldsymbol{x},\boldsymbol{y}) \in \boldsymbol{\mathcal{H}}^2, \quad \big\langle \boldsymbol{\Gamma B x} - \boldsymbol{\Gamma B y} \,\big|\, \boldsymbol{x} - \boldsymbol{y} \big\rangle_{\boldsymbol{\Gamma}^{-1}\boldsymbol{W}} \geq \phi(\|\boldsymbol{\Sigma x} - \boldsymbol{\Sigma y}\|) \,.$$

(v) $\boldsymbol{\Gamma B}$ *is* $\big\|L^{1/2}\Gamma L^{1/2}\big\|^{-1}$*-cocoercive in* $\boldsymbol{\mathcal{H}}_{\boldsymbol{\Gamma}^{-1}\boldsymbol{W}}$.

*Proof.* Let $(\boldsymbol{x},\boldsymbol{y}) \in \boldsymbol{\mathcal{H}}^2$, $x \overset{\text{def}}{=} \boldsymbol{\Sigma x} = \sum_{i=1}^{n} W_i x_i$ and $y \overset{\text{def}}{=} \boldsymbol{\Sigma y} = \sum_{i=1}^{n} W_i y_i$.
(i). Immediate with, respectively, (P1) (ii), (P2) (ii) and (P2) (iii) with Remark 2.5.
(ii). For all $i \in \{1,\ldots,n\}$, $P_{\mathcal{H}_i}$ is linear and continuous, hence so is $\Sigma^*$, and $\boldsymbol{\mathcal{S}} = \Sigma^*(\mathcal{H})$ is thus a closed subspace of $\boldsymbol{\mathcal{H}}$, hence of $\boldsymbol{\mathcal{H}}_{\boldsymbol{\Gamma}^{-1}\boldsymbol{W}}$. Now by definition, $\Sigma^*\boldsymbol{\Sigma x} = \Sigma^* x \in \boldsymbol{\mathcal{S}}$, and using successively Lemma 5.3 (iv) (i), we have for all $\boldsymbol{s} \in \boldsymbol{\mathcal{S}}$, $\langle \boldsymbol{s} \,|\, \boldsymbol{x} - \Sigma^*\boldsymbol{\Sigma x} \rangle_{\boldsymbol{\Gamma}^{-1}\boldsymbol{W}} = \langle \Gamma^{-1}\boldsymbol{\Sigma s} \,|\, x - x \rangle = 0$. Altogether, this characterizes $\Sigma^*\boldsymbol{\Sigma x} = P_{\boldsymbol{\mathcal{S}}}^{\boldsymbol{\Gamma}^{-1}\boldsymbol{W}}\boldsymbol{x}$.
(iii). Let moreover $(\boldsymbol{u},\boldsymbol{v}) \in \boldsymbol{\mathcal{H}}^2$. Supposing that $\boldsymbol{x} \in \big(\mathbf{Id} + \boldsymbol{W}^{-1}\boldsymbol{\Gamma A}\big)^{-1}\boldsymbol{u}$, then by definition $\boldsymbol{u} - \boldsymbol{x} \in \boldsymbol{W}^{-1}\boldsymbol{\Gamma A x}$, that is to say for all $i \in \{1,\ldots,n\}$, $u_i - x_i \in W_i^{-1}\Gamma A_i x_i$. We deduce that

$$\boldsymbol{x} \in J_{\boldsymbol{W}^{-1}\boldsymbol{\Gamma A}}\boldsymbol{u} \Leftrightarrow \forall\, i \in \{1,\ldots,n\},\ x_i \in J_{W_i^{-1}\Gamma A_i} u_i \,. \tag{5.1}$$

Now, for all $i \in \{1,\ldots,n\}$, recall from Remark 2.1 that $A_{i|\mathcal{H}_i}$ is maximally monotone in $\mathcal{H}_i$ and from Remark 2.5 that $\Gamma^{-1}W_i \in \mathfrak{S}_{++}(\mathcal{H}_i)$; we deduce with Combettes and Vũ (2014, Lemma 3.7 (i)) that $W_i^{-1}\Gamma A_{i|\mathcal{H}_i}$ is maximally monotone in $\mathcal{H}_{i\,\Gamma^{-1}W_i}$. By Lemma 5.2 (i)⇔(iv), $J_{W_i^{-1}\Gamma A_i|\mathcal{H}_i}$ has full domain, and with (5.1), so has $J_{\boldsymbol{W}^{-1}\boldsymbol{\Gamma A}}$. Now, supposing moreover $\boldsymbol{y} \in J_{\boldsymbol{W}^{-1}\boldsymbol{\Gamma A}}\boldsymbol{v}$, calling on (5.1), and summing Lemma 5.2 (iii)⇔(iv) over $i \in \{1,\ldots,n\}$, we get $\sum_{i=1}^{n} \langle u_i - v_i \,|\, x_i - y_i \rangle_{\Gamma^{-1}W_i} \geq \sum_{i=1}^{n} \|x_i - y_i\|^2_{\Gamma^{-1}W_i}$, thus $\langle \boldsymbol{u} - \boldsymbol{v} \,|\, \boldsymbol{x} - \boldsymbol{y} \rangle_{\boldsymbol{\Gamma}^{-1}\boldsymbol{W}} \geq \|\boldsymbol{x} - \boldsymbol{y}\|^2_{\boldsymbol{\Gamma}^{-1}\boldsymbol{W}}$. Again with Lemma 5.2, $J_{\boldsymbol{W}^{-1}\boldsymbol{\Gamma A}}$ is also firmly nonexpansive. It is thus the resolvent of a maximally monotone operator in $\boldsymbol{\mathcal{H}}_{\boldsymbol{\Gamma}^{-1}\boldsymbol{W}}$; but the latter is uniquely defined by $J_{\boldsymbol{W}^{-1}\boldsymbol{\Gamma A}}{}^{-1} - \mathbf{Id} = \boldsymbol{W}^{-1}\boldsymbol{\Gamma A}$.
(iv). Since $\operatorname{ran}\boldsymbol{B} \subseteq \boldsymbol{\mathcal{S}}$, developing $\boldsymbol{B}$ and using successively Lemma 5.3 (iii) (iv) (ii) leads to $\big\langle \boldsymbol{\Gamma B x} - \boldsymbol{\Gamma B y} \,\big|\, \boldsymbol{x} - \boldsymbol{y} \big\rangle_{\boldsymbol{\Gamma}^{-1}\boldsymbol{W}} = \big\langle \Gamma^{-1}\boldsymbol{\Sigma}\Sigma^*\Gamma(B\boldsymbol{\Sigma x} - B\boldsymbol{\Sigma y}) \,\big|\, \boldsymbol{\Sigma x} - \boldsymbol{\Sigma y} \big\rangle = \big\langle Bx - By \,\big|\, x - y \big\rangle$. Thus, if $B$ is uniformly monotone with modulus $\phi$, we get $\big\langle \boldsymbol{\Gamma B x} - \boldsymbol{\Gamma B y} \,\big|\, \boldsymbol{x} - \boldsymbol{y} \big\rangle_{\boldsymbol{\Gamma}^{-1}\boldsymbol{W}} \geq \phi(\|x - y\|)$.
(v). Let $\beta \in \mathbb{R}_{++}$. First, using assumption (H1) in the inner product developed above, we get $\big\langle \beta\boldsymbol{\Gamma B x} - \beta\boldsymbol{\Gamma B y} \,\big|\, \boldsymbol{x} - \boldsymbol{y} \big\rangle_{\boldsymbol{\Gamma}^{-1}\boldsymbol{W}} \geq \beta\langle L^{-1}(Bx - By) \,|\, Bx - By\rangle$. Developing again the inner product as above, we also get $\big\|\beta\boldsymbol{\Gamma B x} - \beta\boldsymbol{\Gamma B y}\big\|^2_{\boldsymbol{\Gamma}^{-1}\boldsymbol{W}} = \beta^2\big\langle Bx - By \,\big|\, \Gamma(Bx - By)\big\rangle$. Altogether, denoting $z \overset{\text{def}}{=} L^{-1/2}(Bx - By)$, we obtain

$$\big\langle \beta\boldsymbol{\Gamma B x} - \beta\boldsymbol{\Gamma B y} \,\big|\, \boldsymbol{x} - \boldsymbol{y} \big\rangle_{\boldsymbol{\Gamma}^{-1}\boldsymbol{W}} - \big\|\beta\boldsymbol{\Gamma B x} - \beta\boldsymbol{\Gamma B y}\big\|^2_{\boldsymbol{\Gamma}^{-1}\boldsymbol{W}} \geq \beta\|z\|^2 - \beta^2\big\langle z \,\big|\, L^{1/2}\Gamma L^{1/2} z\big\rangle \,.$$

By Cauchy-Schwartz inequality, $\big\langle z \,\big|\, L^{1/2}\Gamma L^{1/2} z\big\rangle \leq \big\|L^{1/2}\Gamma L^{1/2}\big\|\|z\|^2$, so that with $\beta \overset{\text{set}}{=} \big\|L^{1/2}\Gamma L^{1/2}\big\|^{-1}$, we finally obtain $\big\langle \beta\boldsymbol{\Gamma B x} - \beta\boldsymbol{\Gamma B y} \,\big|\, \boldsymbol{x} - \boldsymbol{y} \big\rangle_{\boldsymbol{\Gamma}^{-1}\boldsymbol{W}} \geq \big\|\beta\boldsymbol{\Gamma B x} - \beta\boldsymbol{\Gamma B y}\big\|^2_{\boldsymbol{\Gamma}^{-1}\boldsymbol{W}}$, and Lemma 5.2 (i)⇔(iii) terminates the proof. ■

## 5.3 Fixed-Point Equation

Now that we have all necessary notions, let us characterize the solutions of (1.1) with the following fixed-point equation.

**Proposition 5.2.** *Under assumptions (H2) and (P1)–(P2), $x \in \mathcal{H}$ is a solution of (1.1) if, and only if, there exists $(z_i)_{1\leq i\leq n} \in \bigtimes_{i=1}^n \mathcal{H}_i$ such that $x = \sum_{i=1}^n W_i z_i$, and for all $i \in \{1,\dots,n\}$*

$$z_i = R_{W_i^{-1}\Gamma A_i}\big((2x - \Gamma Bx)^{\mathcal{H}_i} - z_i\big) - (\Gamma Bx)^{\mathcal{H}_i} . \quad (5.2)$$

*Proof.* Let $x \in \mathcal{H}$. We have the equivalence

$$0 \in \sum_{i=1}^n A_i x + Bx \Leftrightarrow \exists (z_i)_{1\leq i\leq n} \in \bigtimes_{i=1}^n \mathcal{H}_i, \begin{cases} \forall\, i \in \{1,\dots,n\},\ W_i(x - \Gamma Bx - z_i) \in \Gamma A_i x , \\ \text{and } x = \sum_{i=1}^n W_i z_i . \end{cases} \quad (5.3)$$

Indeed, suppose that for all $i \in \{1,\dots,n\}$, $y_i \in A_i x$, such that $0 = \sum_{i=1}^n y_i + Bx$. For all $i \in \{1,\dots,n\}$, (H2) shows that $y_i \in \mathcal{H}_i$, and thanks to (P1) (ii) and (P2) (ii) we can define $z_i \overset{\text{def}}{=} (x - \Gamma Bx)^{\mathcal{H}_i} - W_i^{-1}\Gamma y_i \in \mathcal{H}_i$. Using (P2) (i), we get for all $i \in \{1,\dots,n\}$, $W_i(x - \Gamma Bx - z_i) = W_i\big((x - \Gamma Bx)^{\mathcal{H}_i^\perp} + W_i^{-1}\Gamma y_i\big) = \Gamma y_i \in \Gamma A_i x$, and using in addition (P2) (iv), we get $\sum_{i=1}^n W_i z_i = \sum_{i=1}^n W_i\big((x - \Gamma Bx)^{\mathcal{H}_i} - W_i^{-1}\Gamma y_i\big) = \sum_{i=1}^n W_i(x - \Gamma Bx) - \sum_{i=1}^n \Gamma y_i = x - \Gamma(Bx + \sum_{i=1}^n y_i) = x$. Conversely, summing the right-hand side of (5.3) gives $\sum_{i=1}^n W_i(x - \Gamma Bx - z_i) \in \sum_{i=1}^n \Gamma A_i x$, and using (P2) (iv), we get $x - \Gamma Bx - x \in \sum_{i=1}^n \Gamma A_i x$, and invertibility of $\Gamma$ leads to $0 \in \sum_{i=1}^n A_i x + Bx$. Now thanks to (H2) and (P2) (i), for all $i \in \{1,\dots,n\}$,

$$W_i(x - \Gamma Bx - z_i) \in \Gamma A_i x \Leftrightarrow W_i\big((x - \Gamma Bx)^{\mathcal{H}_i} - z_i\big) \in \Gamma A_i x^{\mathcal{H}_i},$$

using (P1) (ii) and (P2) (ii),

$$\begin{aligned} &\Leftrightarrow (x - \Gamma Bx)^{\mathcal{H}_i} - z_i + x^{\mathcal{H}_i} \in W_i^{-1}\Gamma A_i x^{\mathcal{H}_i} + x^{\mathcal{H}_i} , \\ &\Leftrightarrow x^{\mathcal{H}_i} = J_{W_i^{-1}\Gamma A_i}\big((2x - \Gamma Bx)^{\mathcal{H}_i} - z_i\big) , \\ &\Leftrightarrow z_i = 2J_{W_i^{-1}\Gamma A_i}\big((2x - \Gamma Bx)^{\mathcal{H}_i} - z_i\big) - \big((2x - \Gamma Bx)^{\mathcal{H}_i} - z_i\big) \\ &\qquad\qquad - (\Gamma Bx)^{\mathcal{H}_i} , \\ &\Leftrightarrow z_i = R_{W_i^{-1}\Gamma A_i}\big((2x - \Gamma Bx)^{\mathcal{H}_i} - z_i\big) - (\Gamma Bx)^{\mathcal{H}_i} . \end{aligned}$$

■

The following proposition describes how this translates on the product space $\boldsymbol{\mathcal{H}}$.

**Proposition 5.3.** *Under assumptions (H2) and (P1)–(P2), $x \in \mathcal{H}$ is a solution of (1.1) if, and only if, there exists $\boldsymbol{z} \in \boldsymbol{\mathcal{H}}$ such that $x = \boldsymbol{\Sigma z}$, and $\boldsymbol{z}$ is a fixed point of the operator $\boldsymbol{T}\colon \boldsymbol{\mathcal{H}} \to \boldsymbol{\mathcal{H}}$ defined by*

$$\boldsymbol{T} \overset{\text{def}}{=} \tfrac{1}{2}\big(R_{\boldsymbol{W}^{-1}\boldsymbol{\Gamma A}} R_{\boldsymbol{S}}^{\boldsymbol{\Gamma}^{-1}\boldsymbol{W}} + \mathbf{Id}\big)\big(\mathbf{Id} - \boldsymbol{\Gamma B}\big) .$$

*Proof.* Let $\boldsymbol{z} \overset{\text{def}}{=} (z_i)_{1\leq i\leq n} \in \boldsymbol{\mathcal{H}}$ and $x \overset{\text{def}}{=} \boldsymbol{\Sigma z} = \sum_{i=1}^n W_i z_i$. Since $\boldsymbol{\Gamma B z} = \big((\Gamma Bx)^{\mathcal{H}_i}\big)_{1\leq i\leq n} \in \boldsymbol{\mathcal{S}}$, and with Proposition 5.1 (ii), $P_{\boldsymbol{S}}^{\boldsymbol{\Gamma}^{-1}\boldsymbol{W}} \boldsymbol{z} = (x^{\mathcal{H}_i})_{1\leq i\leq n}$, we obtain $R_{\boldsymbol{S}}^{\boldsymbol{\Gamma}^{-1}\boldsymbol{W}}\big(\mathbf{Id} - \boldsymbol{\Gamma B}\big)\boldsymbol{z} = 2P_{\boldsymbol{S}}^{\boldsymbol{\Gamma}^{-1}\boldsymbol{W}}\boldsymbol{z} - \boldsymbol{z} - \boldsymbol{\Gamma B z} = \big((2x - \Gamma Bx)^{\mathcal{H}_i} - z_i\big)_{1\leq i\leq n}$. Then, Proposition 5.1 (iii) shows that

$$\begin{aligned} \forall\, i \in \{1,\dots,n\},\ z_i \text{ satisfies (5.2)} &\Leftrightarrow \boldsymbol{z} = R_{\boldsymbol{W}^{-1}\boldsymbol{\Gamma A}} R_{\boldsymbol{S}}^{\boldsymbol{\Gamma}^{-1}\boldsymbol{W}}\big(\mathbf{Id} - \boldsymbol{\Gamma B}\big)\boldsymbol{z} - \boldsymbol{\Gamma B z} , \\ &\Leftrightarrow 2\boldsymbol{z} = R_{\boldsymbol{W}^{-1}\boldsymbol{\Gamma A}} R_{\boldsymbol{S}}^{\boldsymbol{\Gamma}^{-1}\boldsymbol{W}}\big(\mathbf{Id} - \boldsymbol{\Gamma B}\big)\boldsymbol{z} + \big(\mathbf{Id} - \boldsymbol{\Gamma B}\big)\boldsymbol{z} , \\ &\Leftrightarrow \boldsymbol{z} = \tfrac{1}{2}\big(R_{\boldsymbol{W}^{-1}\boldsymbol{\Gamma A}} R_{\boldsymbol{S}}^{\boldsymbol{\Gamma}^{-1}\boldsymbol{W}} + \mathbf{Id}\big)\big(\mathbf{Id} - \boldsymbol{\Gamma B}\big)\boldsymbol{z} ; \end{aligned}$$

Proposition 5.2 terminates the proof. ■

The following proposition gives the crucial properties of the above fixed-point operator.

**Proposition 5.4.** *Define the operator* $\boldsymbol{T_A}\colon \mathcal{H} \to \mathcal{H}$ *and* $\boldsymbol{T_B}\colon \mathcal{H} \to \mathcal{H}$ *by*

$$\boldsymbol{T_A} \overset{\text{def}}{=} \tfrac{1}{2}\left(R_{\boldsymbol{W}^{-1}\boldsymbol{\Gamma}\boldsymbol{A}} R_{\boldsymbol{S}}^{\boldsymbol{\Gamma}^{-1}\boldsymbol{W}} + \mathbf{Id}\right) \qquad \text{and} \qquad \boldsymbol{T_B} \overset{\text{def}}{=} \left(\mathbf{Id} - \boldsymbol{\Gamma}\boldsymbol{B}\right).$$

*Under assumptions (H1)–(H2) and (P1)–(P2),*

(i) $\boldsymbol{T_A} \in \mathcal{A}\left(\tfrac{1}{2}, \mathcal{H}_{\boldsymbol{\Gamma}^{-1}\boldsymbol{W}}\right)$,

(ii) $\boldsymbol{T_B} \in \mathcal{A}\left(\tfrac{1}{2}\|L^{1/2}\Gamma L^{1/2}\|, \mathcal{H}_{\boldsymbol{\Gamma}^{-1}\boldsymbol{W}}\right)$, *and*

(iii) $\boldsymbol{T} \in \mathcal{A}\left(\left(2 - \tfrac{1}{2}\|L^{1/2}\Gamma L^{1/2}\|\right)^{-1}, \mathcal{H}_{\boldsymbol{\Gamma}^{-1}\boldsymbol{W}}\right)$.

*Proof.* (i). With Proposition 5.1 (iii) and Lemma 5.2 (ii)⇔(iv), $R_{\boldsymbol{W}^{-1}\boldsymbol{\Gamma}\boldsymbol{A}}$ and $R_{\boldsymbol{S}}^{\boldsymbol{\Gamma}^{-1}\boldsymbol{W}}$ are both non-expansive in $\mathcal{H}_{\boldsymbol{\Gamma}^{-1}\boldsymbol{W}}$, and so is their composition. The result follows by definition. (ii). With Proposition 5.1 (v) and thanks to (P1) (i), Lemma 5.1 (iv) can be applied with $\beta \overset{\text{set}}{=} \|L^{1/2}\Gamma L^{1/2}\|^{-1}$ and $\gamma \overset{\text{set}}{=} 1$. (iii). From $\boldsymbol{T} = \boldsymbol{T_A}\boldsymbol{T_B}$, combine (i), (ii) and Lemma 5.1 (iii) with $\alpha \overset{\text{set}}{=} \tfrac{1}{2}$ and $\alpha' \overset{\text{set}}{=} \tfrac{1}{2}\|L^{1/2}\Gamma L^{1/2}\|$. ■

## 5.4 Convergence

We are now ready to prove our main result.

**Theorem 5.1.** *Set* $\boldsymbol{z}_0 \in \mathcal{H}$ *and define* $(\boldsymbol{z}_k)_{k\in\mathbb{N}}$ *the sequence in* $\mathcal{H}$ *such that for all* $k \in \mathbb{N}$,

$$\boldsymbol{z}_{k+1} = \boldsymbol{z}_k + \rho_k\left(\boldsymbol{T_A}\left(\boldsymbol{T_B}\boldsymbol{z}_k + \boldsymbol{b}_k\right) + \boldsymbol{a}_k - \boldsymbol{z}_k\right), \tag{5.4}$$

*where* $(\boldsymbol{a}_k, \boldsymbol{b}_k) \in \mathcal{H}^2$ *and* $\rho_k \in \left]0, 2 - \tfrac{1}{2}\|L^{1/2}\Gamma L^{1/2}\|\right[$.
*Under assumptions (H1)–(H3) and (P1)–(P2), if*

(i) $\sum_{k\in\mathbb{N}} \rho_k\left(2 - \tfrac{1}{2}\|L^{1/2}\Gamma L^{1/2}\| - \rho_k\right) = +\infty$, *and*

(ii) $\sum_{k\in\mathbb{N}} \rho_k\|\boldsymbol{a}_k\| < +\infty$ *and* $\sum_{k\in\mathbb{N}} \rho_k\|\boldsymbol{b}_k\| < +\infty$,

*then there exists* $\boldsymbol{z} \in \operatorname{fix} \boldsymbol{T}$ *such that* $x \overset{\text{def}}{=} \boldsymbol{\Sigma}\boldsymbol{z}$ *is a solution of (1.1) and that*

(a) $(\boldsymbol{z}_k)_{k\in\mathbb{N}}$ *converges weakly to* $\boldsymbol{z}$, *and*

(b) $(\boldsymbol{\Sigma}\boldsymbol{z}_k)_{k\in\mathbb{N}}$ *converges weakly to* $x$.

*If moreover*

(iii) $0 < \inf_{k\in\mathbb{N}} \{\rho_k\} \leq \sup_{k\in\mathbb{N}} \{\rho_k\} \leq 2 - \tfrac{1}{2}\|L^{1/2}\Gamma L^{1/2}\|$, *and*

(iv) $B$ *is uniformly monotone,*

*then*

(c) $(\boldsymbol{\Sigma}\boldsymbol{z}_k)_{k\in\mathbb{N}}$ *converges strongly to* $x$.

*Proof.* For convenience, define $\alpha \overset{\text{def}}{=} (2 - \frac{1}{2}\|L^{1/2}\Gamma L^{1/2}\|)^{-1}$.

(a)–(b). We know from Proposition 5.4 (iii) that there exists $\boldsymbol{R}\colon \boldsymbol{\mathcal{H}} \to \boldsymbol{\mathcal{H}}$ nonexpansive in $\boldsymbol{\mathcal{H}}_{\Gamma^{-1}W}$ such that $\boldsymbol{T} = \boldsymbol{T_A}\boldsymbol{T_B} = \alpha\boldsymbol{R} + (1-\alpha)\mathbf{Id}$. Then, recurrence (5.4) reads for all $k \in \mathbb{N}$ as

$$\begin{aligned} \boldsymbol{z}_{k+1} &= \boldsymbol{z}_k + \rho_k\big(\alpha\boldsymbol{R}\boldsymbol{z}_k + (1-\alpha)\boldsymbol{z}_k + \boldsymbol{c}_k - \boldsymbol{z}_k\big)\ , \\ &= \boldsymbol{z}_k + \sigma_k\big(\boldsymbol{R}\boldsymbol{z}_k + \boldsymbol{d}_k - \boldsymbol{z}_k\big)\ , \end{aligned}$$

where $\boldsymbol{c}_k \overset{\text{def}}{=} \boldsymbol{T_A}\big(\boldsymbol{T_B}\boldsymbol{z}_k + \boldsymbol{b}_k\big) + \boldsymbol{a}_k - \boldsymbol{T_A}\boldsymbol{T_B}\boldsymbol{z}_k$, $\sigma_k \overset{\text{def}}{=} \alpha\rho_k$ and $\boldsymbol{d}_k \overset{\text{def}}{=} \alpha^{-1}\boldsymbol{c}_k$. By the triangular inequality and nonexpansivity of the involved operators, we have for all $k \in \mathbb{N}$,

$$\|\boldsymbol{c}_k\|_{\Gamma^{-1}W} \le \|\boldsymbol{a}_k\|_{\Gamma^{-1}W} + \|\boldsymbol{T_B}\boldsymbol{z}_k + \boldsymbol{b}_k - \boldsymbol{T_B}\boldsymbol{z}_k\|_{\Gamma^{-1}W} \le \|\boldsymbol{a}_k\|_{\Gamma^{-1}W} + \|\boldsymbol{b}_k\|_{\Gamma^{-1}W}\ .$$

Thus, (ii) provides that $\sum_{k\in\mathbb{N}} \sigma_k\|\boldsymbol{d}_k\|_{\Gamma^{-1}W} = \sum_{k\in\mathbb{N}} \rho_k\|\boldsymbol{c}_k\|_{\Gamma^{-1}W} < +\infty$ by norms equivalence. Moreover, for all $k \in \mathbb{N}$, $\sigma_k < 1$ and thanks to (i), $\sum_{k\in\mathbb{N}} \sigma_k(1-\sigma_k) = \alpha^2 \sum_{k\in\mathbb{N}} \rho_k(\alpha^{-1} - \rho_k) = +\infty$. Finally, Proposition 5.3 and assumption (H3) ensures $\operatorname{fix}\boldsymbol{T} \neq \emptyset$; but $\operatorname{fix}\boldsymbol{T} = \operatorname{fix}\boldsymbol{R}$, and the results follows from Combettes (2004, Lemma 5.1) together with Proposition 5.3 and the continuity of $\boldsymbol{\Sigma}$.

(c). Let now $\boldsymbol{z} \in \operatorname{fix}\boldsymbol{T}$ be the weak limit of the sequence $(\boldsymbol{z}_k)_{k\in\mathbb{N}}$. Let $k \in \mathbb{N}$; we have from the recursion (5.4) that

$$\begin{aligned} \boldsymbol{z}_{k+1} - \boldsymbol{z} &= (1-\rho_k)(\boldsymbol{z}_k - \boldsymbol{z}) + \rho_k\big(\boldsymbol{T_A}(\boldsymbol{T_B}\boldsymbol{z}_k + \boldsymbol{b}_k) + \boldsymbol{a}_k - \boldsymbol{z}\big)\ , \\ &= (1-\rho_k)(\boldsymbol{z}_k - \boldsymbol{z}) + \rho_k(\boldsymbol{T}\boldsymbol{z}_k - \boldsymbol{T}\boldsymbol{z}) + \rho_k\boldsymbol{c}_k\ , \end{aligned}$$

and we can develop

$$\begin{aligned} \|\boldsymbol{z}_{k+1} - \boldsymbol{z} - \rho_k\boldsymbol{c}_k\|^2_{\Gamma^{-1}W} &= (1-\rho_k)^2\|\boldsymbol{z}_k - \boldsymbol{z}\|^2_{\Gamma^{-1}W} + {\rho_k}^2\|\boldsymbol{T}\boldsymbol{z}_k - \boldsymbol{T}\boldsymbol{z}\|^2_{\Gamma^{-1}W} \\ &\quad + 2(1-\rho_k)\rho_k\langle \boldsymbol{z}_k - \boldsymbol{z} \,|\, \boldsymbol{T}\boldsymbol{z}_k - \boldsymbol{T}\boldsymbol{z}\rangle_{\Gamma^{-1}W}\ . \end{aligned} \tag{5.5}$$

At this point, observe that by nonexpansivity of $\boldsymbol{T_A}$, Proposition 5.4 (iii) and Lemma 5.1 (i),

$$\|\boldsymbol{T}\boldsymbol{z}_k - \boldsymbol{T}\boldsymbol{z}\|^2_{\Gamma^{-1}W} \le \|\boldsymbol{T_B}\boldsymbol{z}_k - \boldsymbol{T_B}\boldsymbol{z}\|^2_{\Gamma^{-1}W} \le \|\boldsymbol{z}_k - \boldsymbol{z}\|^2_{\Gamma^{-1}W} - \theta\|\boldsymbol{\Gamma}\boldsymbol{B}\boldsymbol{z}_k - \boldsymbol{\Gamma}\boldsymbol{B}\boldsymbol{z}\|^2_{\Gamma^{-1}W}\ , \tag{5.6}$$

where $\theta \overset{\text{def}}{=} \frac{2-\|L^{1/2}\Gamma L^{1/2}\|}{\|L^{1/2}\Gamma L^{1/2}\|} > 0$. We proceed now by case analysis. Supposing first that $\rho_k \le 1$, we have $(1-\rho_k)\rho_k \ge 0$. Thus by Cauchy-Schwartz inequality and nonexpansivity of $\boldsymbol{T}$, (5.5) becomes

$$\begin{aligned} \|\boldsymbol{z}_{k+1} - \boldsymbol{z} - \rho_k\boldsymbol{c}_k\|^2_{\Gamma^{-1}W} &\le (1-\rho_k)^2\|\boldsymbol{z}_k - \boldsymbol{z}\|^2_{\Gamma^{-1}W} + {\rho_k}^2\|\boldsymbol{T}\boldsymbol{z}_k - \boldsymbol{T}\boldsymbol{z}\|^2_{\Gamma^{-1}W} \\ &\quad + 2(1-\rho_k)\rho_k\|\boldsymbol{z}_k - \boldsymbol{z}\|^2_{\Gamma^{-1}W}\ , \\ &\le \|\boldsymbol{z}_k - \boldsymbol{z}\|^2_{\Gamma^{-1}W} - {\rho_k}^2\theta\|\boldsymbol{\Gamma}\boldsymbol{B}\boldsymbol{z}_k - \boldsymbol{\Gamma}\boldsymbol{B}\boldsymbol{z}\|^2_{\Gamma^{-1}W}\ , \end{aligned} \tag{5.7}$$

where we used (5.6) and $(1-\rho_k)^2 + {\rho_k}^2 + 2(1-\rho_k)\rho_k = (1-\rho_k+\rho_k)^2 = 1$.

Supposing now that $\rho_k > 1$, we have $(1-\rho_k)\rho_k < 0$. We use then Lemma 5.1 (ii) in (5.5) to get

$$\begin{aligned} \|\boldsymbol{z}_{k+1} - \boldsymbol{z} - \rho_k\boldsymbol{c}_k\|^2_{\Gamma^{-1}W} &\le (1-\rho_k)^2\|\boldsymbol{z}_k - \boldsymbol{z}\|^2_{\Gamma^{-1}W} + {\rho_k}^2\|\boldsymbol{T}\boldsymbol{z}_k - \boldsymbol{T}\boldsymbol{z}\|^2_{\Gamma^{-1}W} \\ &\quad + \frac{(1-\rho_k)\rho_k}{1-\alpha}\Big(\|\boldsymbol{T}\boldsymbol{z}_k - \boldsymbol{T}\boldsymbol{z}\|^2_{\Gamma^{-1}W} + (1-2\alpha)\|\boldsymbol{z}_k - \boldsymbol{z}\|^2_{\Gamma^{-1}W}\Big)\ . \end{aligned}$$

In the right-hand side, the factor in front of $\|\boldsymbol{T}\boldsymbol{z}_k - \boldsymbol{T}\boldsymbol{z}\|^2_{\Gamma^{-1}W}$ is ${\rho_k}^2 + \frac{(1-\rho_k)\rho_k}{1-\alpha} = \rho_k\frac{1-\sigma_k}{1-\alpha}$. Recalling that $\sigma_k < 1$, this factor is strictly positive, and we can call on (5.6) to get

$$\|\boldsymbol{z}_{k+1} - \boldsymbol{z} - \rho_k\boldsymbol{c}_k\|^2_{\Gamma^{-1}W} \le \|\boldsymbol{z}_k - \boldsymbol{z}\|^2_{\Gamma^{-1}W} - \frac{1-\sigma_k}{1-\alpha}\rho_k\theta\|\boldsymbol{\Gamma}\boldsymbol{B}\boldsymbol{z}_k - \boldsymbol{\Gamma}\boldsymbol{B}\boldsymbol{z}\|^2_{\Gamma^{-1}W}\ , \tag{5.8}$$

since the factor in front of $\|z_k - z\|^2_{\Gamma^{-1}W}$ simplifies as $(1-\rho_k)^2 + \rho_k{}^2 + \frac{(1-\rho_k)\rho_k}{1-\alpha} + \frac{(1-\rho_k)\rho_k}{1-\alpha}(1-2\alpha) = (1-\rho_k)^2 + \rho_k{}^2 + 2(1-\rho_k)\rho_k = 1$. Altogether, we can combine (5.7) and (5.8), and again Cauchy-Schwartz inequality, to conclude that

$$\tau_k \|\Gamma B z_k - \Gamma B z\|^2_{\Gamma^{-1}W} \leq \|z_k - z\|^2_{\Gamma^{-1}W} - \|z_{k+1} - z\|^2_{\Gamma^{-1}W} - \rho_k{}^2 \|c_k\|^2_{\Gamma^{-1}W} + 2\rho_k \|c_k\|_{\Gamma^{-1}W} \|z_{k+1} - z\|_{\Gamma^{-1}W} \,.$$

where $\tau_k \overset{\text{def}}{=} \min\left\{\frac{1-\sigma_k}{1-\alpha}, \rho_k\right\}\rho_k\theta$. Recall that $(z_k)_{k\in\mathbb{N}}$ converges weakly to $z$ and is therefore bounded (Riesz and Sz.-Nagy, 1990, § 31). Also, $(\rho_k c_k)_{n\in\mathbb{N}}$ is summable by (ii), and $\inf_{k\in\mathbb{N}}\{\tau_k\} > 0$ by (iii). We finally obtain that $\sum_{k\in\mathbb{N}} \|\Gamma B z_k - \Gamma B z\|^2_{\Gamma^{-1}W} < +\infty$, and in particular $(\Gamma B z_k)_{k\in\mathbb{N}}$ converges strongly to $\Gamma B z$, hence $\lim_{k\to\infty} \langle \Gamma B z_k - \Gamma B z \,|\, z_k - z\rangle_{\Gamma^{-1}W} = 0$. If now $B$ is uniformly monotone with modulus $\phi$, it follows from Proposition 5.1 (iv) that $\lim_{k\to\infty} \phi(\|\Sigma z_k - \Sigma z\|) = 0$, and the result follows from the property of $\phi$. ■

**Remark 5.1.** Since $T_A \in \mathcal{A}\left(\frac{1}{2}, \mathcal{H}_{\Gamma^{-1}W}\right)$, there exists an operator $A' \colon \mathcal{H} \to 2^{\mathcal{H}}$ maximally monotone in $\mathcal{H}_{\Gamma^{-1}W}$ such that $T_A = J_{A'}$. We can then show that $\operatorname{fix} T = \operatorname{zer}\{A' + \Gamma B\}$, so that our algorithm actually reduces to a forward-backward splitting on the product space. Following this idea, we discussed in a previous work other conditions of strong convergence (Raguet *et al.*, 2013, Remark 4.2) derived from known results for the forward-backward splitting. Moreover, it is easy to use this formulation to write an *inertial* version of our preconditioned generalized forward-backward algorithm from the work of Lorenz and Pock (2015). However, the resulting algorithm requires some more computations and more than twice as much memory load in comparison to Algorithm 1, whereas numerical experiments showed no improvement in the convergence speed (data not shown).

At last, we can now derive the convergence conditions given in § 2.3.

**Corollary 5.1.** *Theorem 2.1 and Corollary 2.1 hold.*

*Proof.* Skipping some calculations, recursion (2.2) is a specific instance of (5.4), leading to the first claim. Suppose now that assumptions (h1)–(h3) are in force, and set $B \overset{\text{set}}{=} \nabla f$ and for all $i \in \{1, \ldots, n\}$, $A_i \overset{\text{set}}{=} \partial g_i$.
(h1)⇒(H1). The functional $f \circ L^{-1/2}$ is also convex and everywhere differentiable, and for all $(x,y) \in \mathcal{H}^2$, $\langle L^{-1/2}\nabla f L^{-1/2}x \,|\, y - x\rangle = \langle \nabla f L^{-1/2}x \,|\, L^{-1/2}y - L^{-1/2}x\rangle \geq f(L^{-1/2}y) - f(L^{-1/2}x)$, showing that $\nabla(f \circ L^{-1/2}) = L^{-1/2}\nabla f L^{-1/2}$. According to Remark 2.2, $\nabla(f \circ L^{-1/2})$ is 1-Lipschitz continuous in $\mathcal{H}$, hence by the theorem of Baillon and Haddad (1977, Corollaire 10) it is 1-cocoercive in $\mathcal{H}$. For $(x,y) \in \mathcal{H}^2$, set $u \overset{\text{def}}{=} L^{1/2}x$, $v \overset{\text{def}}{=} L^{1/2}y$, and use Lemma 5.2 (i)⇔(iii) to get

$$\langle \nabla(f \circ L^{-1/2})u - \nabla(f \circ L^{-1/2})v \,|\, u - v\rangle^2 \geq \|\nabla(f \circ L^{-1/2})u - \nabla(f \circ L^{-1/2})v\|^2 \,,$$

so that

$$\begin{aligned}\langle \nabla f x - \nabla f y \,|\, x - y\rangle^2 &\geq \|L^{-1/2}\nabla f x - L^{-1/2}\nabla f y\|^2 \,,\\ &\geq \|\nabla f x - \nabla f y\|^2_{L^{-1}} \,.\end{aligned}$$

(h2)⇒(H2). Let $i \in \{1, \ldots, n\}$. Moreau (1965, Proposition 12.b) shows that $\partial g_i$ is maximally monotone in $\mathcal{H}$. Moreover, since $g_i = g_i \circ P_{\mathcal{H}_i}$ is proper, we know that $\operatorname{dom} g_i \cap \mathcal{H}_i \neq \emptyset$; since $\operatorname{ran} P_{\mathcal{H}_i} = \mathcal{H}_i$ is a closed subspace, Bauschke and Combettes (2011, Corollary 16.37) show that $\partial g_i = \partial(g_i \circ P_{\mathcal{H}_i}) = P_{\mathcal{H}_i} \partial g_i P_{\mathcal{H}_i}$.

(h3)⇒(H3). (h3) (i) together with Bauschke and Combettes (2011, Corollary 16.37) show that $\partial(f + \sum_{i=1}^n g_i) = \nabla f + \sum_{i=1}^n \partial g_i$. But $\operatorname{zer}\{\partial(f + \sum_{i=1}^n g_i)\} = \operatorname{argmin}\{f + \sum_{i=1}^n g_i\}$, according to Fermat principle.
To obtain the final formulation of Corollary 2.1, we follow Hiriart-Urruty and Lemaréchal (1993, § XV.4) with Remark 2.1 to get for all $i \in \{1, \ldots, n\}$, $(J_{W_i^{-1}\Gamma\partial g_i})_{|\mathcal{H}_i} = J_{W_i^{-1}\Gamma\partial g_{i|\mathcal{H}_i}} = \operatorname{prox}^{\Gamma^{-1}W_i}_{g_{i|\mathcal{H}_i}}$. At last for the strong convergence, consider that uniform convexity of $f$ implies uniform monotonicity of $\nabla f$ (see Bauschke and Combettes, 2011, Example 22.3 (iii)). ■